\documentclass[11pt]{amsart}
\usepackage[headings]{fullpage}
\usepackage{amsmath}
\usepackage{amsfonts}
\usepackage{latexsym}
\usepackage{amssymb}
\usepackage{MnSymbol}
\newtheorem{thm}{Theorem}[section]

\newtheorem{lemma}{Lemma}[section]
\newtheorem{corr}{Corollary}[section]
\newtheorem{defn}{Definition}[section]

\newtheorem{proposition}{Proposition}[section]
\def\homo{\mathop{\sf Hom}}
\def\ker{\mathop{\sf Ker}}

\def\curl{\mathop{\sf curl}}
\def\rk{\mathop{\sf rank}}
\def\det{\mathop{\sf det}}
\def\ann{\mathop{\sf ann}}
\def\dep{\mathop{\sf depth}}
\def\dim{\mathop{\sf dim}}
\def\pd{\mathop{\sf pd}}
\def\height{\mathop{\sf height}}
\def\ass{\mathop{\sf Ass}}
\def\supp{\mathop{\sf Supp}}
\def\pt{{\sf PT}}
\def\t{{\sf T}}
\def\tor{\mathop{\sf tor}}
\def\min{{\sf min}}

\title[The Hautus Test]{The Hautus Test and Genericity Results for Controllable and Uncontrollable Behaviors}
  
\author[S.Shankar]{Shiva Shankar}

\begin{document}

\maketitle

\begin{center}

{\small Chennai Mathematical Institute,\\
Plot No. H1, SIPCOT IT Park, Kelambakkam,\\
Siruseri, Chennai (Madras) - 603103, India.}\\

\vspace{1cm}

\emph{In Memoriam: Professor E.G.F.~Thomas (1939 - 2011)}

\end{center}

\vspace{.5cm}

\begin{abstract} \noindent The computational effectiveness of Kalman's state space controllability rests on the well known {\it Hautus Test} which describes a rank condition of the matrix $(\frac{d}{dt}I-A, ~B)$. This paper generalizes this test to a generic class of behaviors (belonging to a Zariski open set) defined by systems of PDE ~(i.e. systems which arise as kernels of operators given by matrices $(p_{ij}(\partial))$ whose entries are in $\mathbb{C}[\partial_1, \ldots , \partial_n]$), and studies its implications, especially to issues of genericity. The paper distinguishes two classes of systems, under-determined and over-determined. The Hautus test developed here implies that a generic strictly under-determined system is controllable, whereas a generic over-determined system is uncontrollable.
\end{abstract}

{\tiny \hspace{1cm}AMS classification: 93B05, 35Q93, 13P25} \\

\vspace{1cm}

\section{Introduction} The developments in `post war' Control Theory rest largely on Kalman's innovative translation of the Control Problem to a problem about differential equations which describes the evolution of the state of a linear input-output system. The form of this equation is 

\begin{equation}
\frac{dx}{dt} = Ax + Bu
\end{equation}
where the state $x$ is in $\mathbb{R}^\ell$, the input $u$ is in $\mathbb{R}^m$, and $A:\mathbb{R}^\ell \rightarrow \mathbb{R}^\ell$, ~$B:\mathbb{R}^m \rightarrow \mathbb{R}^\ell$ are linear maps. The system (1) is said to be {\it controllable} (or state controllable to distinguish it from other related notions) if given two states $x_1$ and $x_2$, there is an input $u:I \rightarrow \mathbb{R}^m$ ($I \subset \mathbb{R}$, an interval in time) such that $x(t_1) = x_1, ~x(t_2) =  x_2$, ~for some $t_1, t_2$ in $I$ \cite{k}.     

Kalman's definition has been generalized in many directions, in particular to (affine) nonlinear systems defined on smooth manifolds, and to infinite dimensional systems described by 1-parameter semigroups. This paper is concerned with a more recent and far reaching generalization due to J.C.~Willems, summarized in \cite{w}. (For a historical perspective of this evolution of ideas, please see \cite{se}.) Willems' Theory of Behaviors is designed to overcome foundational problems with the Kalman Theory, such as an {\it a priori} assumption of a causal structure and its attendant division of signals into inputs and outputs, the emphasis on state (i.e. first order equations), and the loss of computational effectiveness in other generalizations of Kalman's finite dimensional linear state space theory.
Here, computational effectiveness of a theory refers to the existence of  finite procedures to verify its statements. For instance, computational effectiveness of state controllability is guaranteed by the well known Hautus test:

\begin{thm} \cite{h} Define $H(\frac{d}{dt}) = (\frac{d}{dt}I_\ell - A, ~ -B)$. Then (1) is (state) controllable if and only if for every $\lambda$ in $\mathbb{C}$, $\rk (H(\lambda))$ is constant (and thus if and only if $H(\lambda)$ has full row rank for every $\lambda $ in $\mathbb{C}$).
\hspace*{\fill}$\square$
\end{thm}

\noindent Remark 1.1: The Hautus test is usually written as the above rank condition for the matrix $(\frac{d}{dt}I_\ell - A, ~B)$. As the two ranks are identical, the above statement is dictated only by notational reasons explained in the next section. \\

Thus, to determine whether (1) is controllable, it suffices to calculate all the maximal minors of $H(\frac{d}{dt})$, i.e. the determinants of the $\ell \times \ell$ submatrices of $H(\frac{d}{dt})$ -- there are ${\tiny \left(\begin{array}{c}\ell+m \\\ell\end{array}\right)}$ of them -- and then to check that these determinants do not have a common zero in $\mathbb{C}$. In other words, is the $\ell$-th determinantal ideal $\frak{i}_\ell$, generated by the maximal minors of $H(\frac{d}{dt})$, equal to $\mathbb{C}[\frac{d}{dt}]$? This computational effectiveness of the Kalman theory was a major reason for its success, and one of the advantages of Willems' behavioral theory is that it retains it in a more general setting. Indeed, this paper shows that behavioral controllability can be effectively determined by a test which is a generalization of the above Hautus test for state space systems. 

In the more general setting of Willems' theory, a dynamical system is identified with the collection of all the trajectories (or signals) that can possibly occur. Controllability in this setting is not the ability to move from one state to another in finite time as in the Kalman theory - indeed there is no notion of state now - but is, instead, the ability to move from one trajectory to another in finite time. In its generalization to behaviors defined by systems of PDE, controllability assumes the following form:  

\begin{defn} Let $\mathcal{A} = \mathbb{C}[\partial_1, \ldots ,\partial_n]$ be the ring of (constant coefficient) partial differential operators, and let 
$\mathcal{P}$ be an $\mathcal{A}$-submodule of $\mathcal{A}^k$. Let 
$P(\partial):\mathcal{F}^k \rightarrow \mathcal{F}^\ell$ be the map defined by an $\ell \times k$ matrix $P(\partial)$ whose $\ell$ rows generate $\mathcal{P}$, and where $\mathcal{F}$ is either the space $\mathcal{D}'$ of distributions or the space $\mathcal{C}^\infty$ of smooth functions on $\mathbb{R}^n$. The behavior $\mathcal{B}_\mathcal{F}(P(\partial))$, given by the kernel $\ker_\mathcal{F}(P(\partial))$ of the above map, is said to be {\it controllable} if given two subsets $U_1$ and $U_2$ of $\mathbb{R}^n$ whose closures do not intersect, and two elements $f_1$ and $f_2$ of the behavior, there is an element $f$ in the behavior such that $f=f_1$ on some neighborhood of $U_1$, and $f=f_2$ on some neighborhood of $U_2$.
\end{defn} 

\noindent Remark 1.2: The behavior $\mathcal{B}_\mathcal{F}(P(\partial))$ in the above definition depends only on the submodule $\mathcal{P}$ -- it is isomorphic to $\homo_\mathcal{A}(\mathcal{A}^k/\mathcal{P},~\mathcal{F})$ -- and will therefore be denoted $\mathcal{B}_\mathcal{F}(\mathcal{P})$. 

\begin{thm} \cite{ps} The behavior $\mathcal{B}_\mathcal{F}(\mathcal{P})$, when $\mathcal{F}$ is either $\mathcal{D}'$ or $\mathcal{C}^\infty $, is controllable if and only if $\mathcal{A}^k/\mathcal{P}$ is torsion free.
\hspace*{\fill}$\square$
\end{thm}

\noindent Remark 1.3: The above theorem is valid when $\mathcal{F}$ is either $\mathcal{D}'$ or $\mathcal{C}^\infty$ because they are {\it injective cogenerators} as $\mathcal{A}$-modules. These are the spaces of principal concern in this paper. Later, other spaces are briefly considered, such as the space $\mathcal{S}'$ of temperate distributions and spaces of periodic functions, where necessary and sufficient conditions for controllability are different (Theorem 3.1 in \cite{sg} and Theorem 4.2 in \cite{nps}, quoted in Section 5 below). All these conditions are actually necessary and sufficient conditions that a behavior, given as the kernel of a map $P(\partial):\mathcal{F}^k \rightarrow \mathcal{F}^\ell$, admit an image representation, i.e. also be equal to the image of some map $M(\partial): \mathcal{F}^{k_1} \rightarrow \mathcal{F}^k$. It is elementary that an image is controllable in the sense of Definition 1.1 \cite{ps}. 

Every behavior in the space $\mathcal{D}$ of compactly supported smooth functions, the space $\mathcal{E}'$ of compactly supported distributions, or the Schwartz space $\mathcal{S}$ of rapidly decreasing functions, admits an image representation (these spaces are {\it flat} $\mathcal{A}$-modules), and is therefore controllable \cite{sg}. Hence a test for controllability is vacuous here.\\

The purpose of this paper is to use a generalization of the Hautus test to study genericity questions about controllablity (or lack of it) of distributed behaviors, i.e. when the ring of differential operators is $\mathbb{C}[\partial_1, \ldots ,\partial_n]$. This generalization is valid for behaviors defined by free submodules $\mathcal{P}$, and counterexamples show that this is the full extent of a Hautus kind of test. When $n = 1$, i.e. $\mathcal{A} = \mathbb{C}[\frac{d}{dt}]$, the ring is a PID, hence every submodule of $\mathcal{A}^k$ being torsion free is free, and so the Hautus test is valid for every lumped behavior. When $n = 2$, i.e. $\mathcal{A} = \mathbb{C}[\partial_1, \partial_2]$, every controllable behavior is given by a free submodule (Corollary 4 in \cite{ps}) and so the Hautus test is again valid for every $2-D$ behavior (Section 3 below). For general $n$, free submodules are {\it generic} (precise statements occur in Section 6), so that the Hautus test developed here is valid for behaviors outside a vanishingly thin set. It turns out that in this case, the rank of the defining matrix $P(\partial)$ can drop, but only at points $\lambda$ in $\mathbb{C}^n$ that lie on an algebraic variety of dimension $n-2$ or smaller. As such sets in $\mathbb{C}$ are empty, this generalization does indeed coincide with the state space Hautus test of Theorem 1.1 above.

This paper distinguishes two classes of behaviors, the under-determined ones and the over-determined. The Hautus test of this paper implies that, {\it a generic strictly under-determined behavior is controllable, whereas a generic over-determined behavior is uncontrollable.} These results could be considered the principal contribution of this paper. 

Hautus tests for delay-differential systems described by operators of the kind $P(\frac{d}{dt}, \Delta)$ ($\Delta$ is the unit delay) have been obtained by Gl\"using-L\"uerssen \cite{g} and by Rocha and Willems \cite{rj}.
There are also other approaches to generalizing the classical Hautus test, for instance Bourles and Marinescu \cite{bm} and Lomadze \cite{l1}, who use ideas from homological algebra.   

As observed above, the torsion-free condition of Theorem 1.2 is also the condition for the behavior, given as a kernel, to admit an image representation (cohomology vanishing), Oberst \cite{o}. (There are other points of contact with Oberst's seminal paper which are highlighted as remarks in this paper.) In Physics, the existence of an image representation is precisely the existence of a potential, Pommaret and Quadrat \cite{pq}. Recently, Lomadze \cite{l} has shown that controllability implies that every system trajectory can be obtained from a `transfer' trajectory by differentiation, where a transfer trajectory is obtained formally from a transfer function. Kalman's notion of controllability, and its generalization by Willems to the setting of behaviors, thus remains an important issue in control theory. \\

\noindent Remark 1.3: A couple of results in this paper admit proofs shorter than those given here, but which rely on some facts about Cohen-Macaulay rings. I have however chosen to give these longer proofs for two reasons: first, the proof itself is used later in the paper; and second, the hope that these elementary methods will make the paper accessible to a larger number of control theorists who may not have the background in commutative algebra that they would otherwise need. I do, however, give references to the shorter proofs in the literature at apposite points in the text of the paper. 

\vspace{.3cm}

\section{The Hautus test for strong controllability} The starting point of this development is to rewrite equation (1) as
\[
\left(\begin{array}{cc}\frac{d}{dt}I_\ell -A, & -B\end{array}\right)\left(\begin{array}{c}x \\u\end{array}\right) = 0 ~;
\]
in other words, to consider all possible trajectories $(x,u)$ of the state space system as the kernel of the map
\[
\begin{array}{lccc}
\left(\begin{array}{cc}\frac{d}{dt}I_\ell -A, & -B\end{array}\right): & \mathcal{F}^k & \longrightarrow & \mathcal{F}^\ell \\ 
& (x,u) & \mapsto & (\frac{d}{dt}I_\ell -A,~-B)\left(\begin{array}{c}x \\u\end{array}\right)
\end{array}
\]
where $\mathcal{F}$ is either $\mathcal{D}'$ or $\mathcal{C}^\infty $, and $k = \ell + m$. (This now explains the remark on notation after Theorem 1.1) More generally, a lumped behavior is given by the kernel of a map
\begin{equation}
P(\frac{d}{dt}): \mathcal{F}^k \longrightarrow \mathcal{F}^\ell
\end{equation}
where the entries of the $\ell \times k$ matrix $P(\frac{d}{dt})$ are from the ring $\mathcal{A} = \mathbb{C}[\frac{d}{dt}]$, and $\ell$ in general could be larger than $k$ (thus, a state space system is a behavior of a very special kind). Let $\mathcal{P}$ be the submodule of $\mathcal{A}^k$ generated by the rows of $P(\frac{d}{dt})$, so that the above behavior, $\mathcal{B}_\mathcal{F}(\mathcal{P})$, is isomorphic to $\homo_\mathcal{A}(\mathcal{A}^k/\mathcal{P}, ~\mathcal{F})$ (Remark 1.2). As the ring $\mathbb{C}[\frac{d}{dt}]$ is a PID, the finitely generated submodule $\mathcal{P}$ being torsion free, is free. If the $\ell$ rows of $P(\frac{d}{dt})$ is a {\it minimum} set of generators for $\mathcal{P}$, minimum in the sense that it cannot be generated by any set of $\ell -1$ elements, then the rows are actually a basis for $\mathcal{P}$, and hence, necessarily, $\ell \leqslant k$. In this case the matrix $P(\frac{d}{dt})$ has full row rank over the field $\mathbb{C}(\frac{d}{dt})$ of `rational symbols'. 
(The 0 submodule which is generated by the empty set of elements, and whose behavior is all of $\mathcal{F}^k$, is excluded from all discussion, as all the considerations of this paper are trivial in this case.)

 The Hautus test generalizes perfectly to arbitrary lumped behaviors: 

\begin{proposition} (Willems \cite{w}) Let $\mathcal{P}$ be a submodule of $\mathcal{A}^k$ and let $P(\frac{d}{dt})$ be any ($\ell \times k$) matrix of full row rank whose rows generate $\mathcal{P}$. Then the behavior $\mathcal{B}_\mathcal{F}(\mathcal{P})$ given by the kernel of (2), $\mathcal{F}$ either $\mathcal{D}'$ or $\mathcal{C}^\infty$, is controllable (in the sense of Definition 1.1 with $n = 1$) if and only if the matrix $P(\frac{d}{dt})$ satisfies the Hautus condition, namely that for every $\lambda $ in $\mathbb{C}$, $P(\lambda)$ has full row rank, or equivalently if and only if $\mathcal{A}^k/\mathcal{P}$ is free. 
\end{proposition}

\noindent Proof: The structure theory for modules over a PID implies that there is a basis $e_1, \ldots , e_k$ for $\mathcal{A}^k$, and elements $a_1, \dots , a_\ell$ of $\mathcal{A}$ such that $a_1e_1, \dots , a_\ell e_\ell$ is a basis for $\mathcal{P}$, namely the Smith form of the matrix $P(\frac{d}{dt})$. Thus, for $\lambda$ in $\mathbb{C}$, the rank of $P(\lambda)$ equals $\ell $ if and only $a_1(\lambda), \dots , a_\ell (\lambda)$ are all nonzero. This latter condition holds for every $\lambda$ in $\mathbb{C}$ if and only if all the $a_j, 1\leqslant j \leqslant \ell$, are constants. In other words, $P(\lambda)$ has full row rank for every $\lambda$ in $\mathbb{C}$ if and only if a basis for $\mathcal{P}$ extends to a basis for $\mathcal{A}^k$.

This last condition is in turn equivalent to saying that the exact sequence
\[
0 \rightarrow \mathcal{P} \longrightarrow \mathcal{A}^k \longrightarrow \mathcal{A}^k/\mathcal{P} \rightarrow 0
\] 
splits (the condition asserts that $ 0 \rightarrow \mathcal{P} \rightarrow \mathcal{A}^k $ splits). This implies that $\mathcal{A}^k/\mathcal{P}$ is a submodule of $\mathcal{A}^k$, hence torsion free, and so free. By 
Theorem 1.2, this is equivalent to the controllability of $\mathcal{B}_\mathcal{F}(\mathcal{P})$.
\hspace*{\fill}$\square$\\

\noindent Remark 2.1: As observed earlier in the introduction, to say that $P(\lambda)$ has full row rank for every $\lambda$ in $\mathbb{C}$ is to say that the $\ell$-th determinantal ideal $\frak{i}_\ell$ of $P(\frac{d}{dt})$  is equal to $\mathbb{C}[\frac{d}{dt}]$. This ideal equals the $(k-\ell)$-th Fitting ideal of $\mathcal{A}^k/\mathcal{P}$ and is therefore independent of the choice of the matrix $P(\frac{d}{dt})$ whose rows generate $\mathcal{P}$.

\begin{corr} (Willems \cite{w}) A state space behavior is controllable in the sense of behaviors if and only if it is state controllable. Thus Willems' definition of behavioral controllability is a faithful generalization of Kalman's state controllability.
\hspace*{\fill}$\square$
\end{corr}

\noindent Remark 2.2: Wolovich considers a system of the form $\{(y,u)|P(\frac{d}{dt})y=Q(\frac{d}{dt})u\}, \det(P) \neq 0$, and shows that it is controllable in the Kalman sense if and only if $P$ and $Q$ are coprime (Theorem 5.3.1 in \cite{wo}). This is equivalent to the Hautus test of Theorem 1.1. \\

Consider now distributed behaviors on $\mathbb{R}^n$, i.e. when the ring $\mathcal{A}$ of differential operators is $\mathbb{C}[\partial_1, \ldots ,\partial_n]$. Proposition 2.1 suggests the following questions about a distributed behavior $\mathcal{B}_\mathcal{F}(\mathcal{P})$, which the rest of the paper addresses: \\

\noindent Let the rows of the matrix $P(\partial)$ be a minimum set of generators for $\mathcal{P}$ (Definition 2.1 below). What does it mean if $P(\lambda)$ has full row rank for each $\lambda$ in $\mathbb{C}^n$? Is this a necessary and sufficient condition for controllability of $\mathcal{B}_\mathcal{F}(\mathcal{P})$? If not, what is the analogue of the Hautus test for distributed behaviors? What is the Hautus test for signal spaces $\mathcal{F}$ other than $\mathcal{D}'$ and $\mathcal{C}^\infty$? What can be said about the set of all controllable behaviors?\\

As $\mathcal{A}$ is now not a PID, a statement analogous to Proposition 2.1 first requires a definition:

\begin{defn} Let $\mathcal{P}$ be a submodule of $\mathcal{A}^k$. Let $\ell$ be the smallest integer such that $\mathcal{P}$ can be generated by some $\ell$ elements. Then any set of generators for $\mathcal{P}$, $\ell$ in number, is said to be minimum. 

Suppose further that $\mathcal{P}$ is a free submodule of $\mathcal{A}^k$. Then any $\ell \times k$ matrix $P(\partial)$, whose rows is a minimum set of generators for $\mathcal{P}$, and which therefore is a basis for $\mathcal{P}$, has full row rank (over the field of rational symbols $\mathbb{C}(\partial_1, \ldots ,\partial_n)$). \end{defn}

The next proposition shows that the Hautus test, translated {\it verbatim} from the case of lumped behaviors in Proposition 2.1, yields a similar result. The proof is a variant of the proof of Proposition 2.1, and will be further generalized below in Theorem 3.1.

\begin{proposition} (Oberst \cite{o}, ~pp. 156-158) Let $\mathcal{P}$ be a submodule of $\mathcal{A}^k$, and let $P(\partial)$ be any $\ell \times k$ matrix whose $\ell$ rows is a minimum set of generators for $\mathcal{P}$. Then $\mathcal{A}^k/\mathcal{P}$ is free if and only if $P(\lambda)$ has full row rank for every $\lambda $ in $\mathbb{C}^n$ (or in other words, the $\ell$-th determinantal ideal of $P(\partial)$ equals $\mathcal{A}$). Thus $\mathcal{B}_\mathcal{F}(\mathcal{P})$, $\mathcal{F}$ either $\mathcal{D}'$ or $\mathcal{C}^\infty$, is controllable if $P(\partial)$ satisfies the Hautus condition. 
\end{proposition}

\noindent Proof: Clearly, either of the above statements implies that $\ell \leqslant k$: for suppose that $\mathcal{A}^k/\mathcal{P}$ is free, then $0 \rightarrow \mathcal{P} \longrightarrow \mathcal{A}^k \longrightarrow \mathcal{A}^k/\mathcal{P} \rightarrow 0$ splits, hence
$\mathcal{P}$ is a direct summand of $\mathcal{A}^k$, and therefore by Quillen-Suslin, free, of rank $\ell \leq k$; the other statement manifestly implies the inequality. Now, the value $p_{ij}(\lambda)$ at $\lambda$ in $\mathbb{C}^n$ of an entry $p_{ij}(\partial)$ of the matrix $P(\partial)$ is its image under the morphism $\mathcal{A} \rightarrow \mathcal{A}/\frak{m}_\lambda$, where $\frak{m}_\lambda$ is the maximal ideal $(\partial _1 - \lambda _1, \ldots , \partial_n - \lambda _n)$ corresponding to $\lambda = (\lambda_1, \ldots , \lambda _n)$. Thus $\rk(P(\lambda)) = \ell$ if and only if there is a minor, say the determinant $\det(M(\lambda))$ of an $\ell \times \ell$ submatrix $M(\lambda)$ of $P(\lambda)$, which is not 0. This implies that $\det(M(\partial))$ does not belong to $\frak{m}_\lambda$. Then in the localization $\mathcal{A}_{\frak{m}_\lambda}$, $M(\partial)$ is invertible, hence
there is a basis such that the matrix of the localization $P_{\frak{m}_\lambda} (\partial):\mathcal{A}^k_{\frak{m}_\lambda} \rightarrow \mathcal{A}^\ell_{\frak{m}_\lambda}$ (of the morphism $P(\partial):\mathcal{A}^k \rightarrow \mathcal{A}^\ell$) has an $\ell \times \ell$ submatrix equal to the identity $I_\ell$. By row and column operations, all other entries of $P_{\frak{m}_\lambda}(\partial)$ can be made zero. This implies that 
\[
0 \rightarrow \mathcal{P}_{\frak{m}_\lambda} \longrightarrow \mathcal{A}^k_{\frak{m}_\lambda} \longrightarrow \mathcal{A}^k_{\frak{m}_\lambda}/\mathcal{P}_{\frak{m}_\lambda} \rightarrow 0
\]
splits. Thus $\mathcal{A}^k_{\frak{m}_\lambda}/\mathcal{P}_{\frak{m}_\lambda}\simeq (\mathcal{A}^k/\mathcal{P})_{\frak{m}_\lambda}$ is projective, and as this is true for every maximal ideal $\frak{m}_\lambda$, it follows that $\mathcal{A}^k/\mathcal{P}$ is projective. Again by Quillen-Suslin, $ \mathcal{A}^k/\mathcal{P}$ is free.
\hspace*{\fill}$\square$\\

By Theorem 1.2, freeness of $\mathcal{A}^k/\mathcal{P}$ is a sufficient but not necessary condition for $\mathcal{B}_\mathcal{F}(\mathcal{P})$ to be controllable. However $\mathcal{A}^k/\mathcal{P}$ is free if and only if it, and therefore also $\mathcal{P}$, are direct summands of $\mathcal{A}^k$. Then $\mathcal{A}^k \simeq \mathcal{P} \oplus \mathcal{A}^k/\mathcal{P}$, and this implies that $\mathcal{B}_\mathcal{F}(\mathcal{P}) \simeq \homo_\mathcal{A}(\mathcal{A}^k/\mathcal{P}, \mathcal{F})$ is a direct summand of $\mathcal{F}^k$. Such a behavior is said to be {\em strongly controllable} and the corresponding $P(\partial)$ {\it zero left prime} (for instance \cite{rw}). Conversely, suppose that $\mathcal{B}_\mathcal{F}(\mathcal{P})$ is a direct summand of $\mathcal{F}^k$, then any complementary summand is also a behavior as $\mathcal{F}$ is injective (for instance \cite{sg}). Let $\mathcal{B}$ be a choice of such a summand. As $\mathcal{F}$ is also a cogenerator (remark after Theorem 1.2), this behavior $\mathcal{B}$ equals $\mathcal{B}_\mathcal{F}(\mathcal{Q})$ for a unique submodule $\mathcal{Q}$ of $\mathcal{A}^k$. It then follows that $\mathcal{A}^k \simeq \mathcal{P} \oplus \mathcal{Q}$ (for instance \cite{sc}). \\

Thus the above proposition gives a sufficient and necessary condition for strong controllability of $\mathcal{B}_\mathcal{F} (\mathcal{P})$. \\

Proposition 2.1 now implies

\begin{corr} Every controllable lumped behavior is strongly controllable.
\hspace*{\fill}$\square$
\end{corr}

\vspace{.22cm}

\section{The Hautus test for controllability of distributed behaviors}
What then are sufficient {\it and} necessary conditions, akin to the Hautus test, for controllability of a distributed behavior? Under generic conditions, made precise later, the next result is such a test.

\begin{thm} Let $\mathcal{P}$ be a submodule of $\mathcal{A}^k$, and let $P(\partial)$ be any $\ell \times k$ matrix whose $\ell$ rows generate $\mathcal{P}$. Suppose that the $\ell$-th determinantal ideal $\frak{i}_\ell$ of $P(\partial)$ is not the zero ideal (so that in particular $\ell \leqslant k$). Then $\mathcal{A}^k/\mathcal{P}$ is torsion free, and hence $\mathcal{B}_\mathcal{F}(\mathcal{P})$ ($\mathcal{F}$ either $\mathcal{D}'$ or $\mathcal{C}^\infty$) is controllable, if and only if $P(\lambda)$ has full row rank for all $\lambda$ in the complement of an algebraic variety in $\mathbb{C}^n$ of dimension $\leqslant n-2$ (or in other words, that the Krull dimension of the ring $\mathcal{A}/\frak{i}_\ell$ be less than or equal to $n-2$).
\end{thm}

\noindent Proof: Suppose to the contrary that $\rk(P(\lambda))$ is less than $\ell$ for $\lambda$ on an algebraic variety of dimension
$n-1$ in $\mathbb{C}^n$ - it cannot be less than $\ell$ on all of $\mathbb{C}^n$ because $\frak{i}_\ell$ is nonzero, by assumption. Let $\mathcal{V}$ be an irreducible component, and let it be the zero locus of the irreducible polynomial $p(\partial)$ (Krull's Principal Ideal Theorem). This means that each of the ${\tiny \left(\begin{array}{c}k \\\ell\end{array}\right)}$ generators of $\frak{i}_\ell$ is divisible by $p(\partial)$. Let $\frak{p}$ be the prime ideal $(p(\partial))$, and let $P_\frak{p}(\partial)$ be the image of the matrix $P(\partial)$ in the localization $\mathcal{A}_\frak{p}$ of $\mathcal{A}$ at $\frak{p}$. Suppose that every entry of some row of $P_\frak{p}(\partial)$, say the first row, is divisible by $p(\partial)$. As divisibility by $p(\partial)$ in $\mathcal{A}$ and $\mathcal{A}_\frak{p}$ are equivalent, the corresponding row of $P(\partial)$ is also divisible by $p(\partial)$, and then clearly $\mathcal{A}^k/\mathcal{P}$ has a torsion element. By Theorem 1.2, $\mathcal{B}_\mathcal{F}(\mathcal{P})$ is then not controllable. Otherwise, at least one element of this first row of $P_\frak{p}(\partial)$, say the first element, is not divisible by $p(\partial)$, and is therefore a unit in $\mathcal{A}_\frak{p}$. By row and column operations, every other element of the first column, and of the first row, can be made zero. Let the matrix so obtained after these row and column operations, with a unit in the (1,1) entry and the other entries of the first row and column equal to 0, be denoted $P^1_\frak{p}(\partial)$, and let $\mathcal{P}^1_{\frak{p}}$ be the submodule of $\mathcal{A}^k_{\frak{p}}$ generated by its rows. 

The above argument for $P_\frak{p}(\partial)$ can be repeated now for $P^1_\frak{p}(\partial)$, and it follows that either $p(\partial)$ divides every element of some row, say the second, or that some element in that row, say  the (2,2) entry (the (2,1) entry is zero) is a unit in $\mathcal{A}_\frak{p}$. In the first case, $\mathcal{A}^k_{\frak{p}}/\mathcal{P}^1_{\frak{p}}$ has a torsion element, hence so does $\mathcal{A}^k/\mathcal{P}$, and $\mathcal{B}_\mathcal{F}(\mathcal{P})$ is not controllable; otherwise all other elements in the second row and second column can be made zero by column and row operations. Eventually, after at most $\ell-1$ such steps, the resultant matrix has units in the $(j,j)$ entries, $1\leqslant j \leqslant \ell-1$, and the other entries zero except in positions $(\ell,j), \ell \leqslant j \leqslant k$. If now $p(\partial)$ does not divide every entry of the $\ell$-th row, then it also does not divide the generators of the determinant ideal $\frak{i}_\ell$, which is a contradiction.

Conversely, suppose $\mathcal{B}_\mathcal{F}(\mathcal{P})$ is not controllable so that $\mathcal{A}^k/\mathcal{P}$ has a torsion element. This implies that there is an element $x(\partial) = (a_1(\partial), \dots ,a_k(\partial))$ in $\mathcal{A}^k \setminus \mathcal{P}$ such that $r(\partial) = p(\partial)x(\partial) = (p(\partial)a_1(\partial), \ldots ,p(\partial)a_k(\partial))$ is in $\mathcal{P}$, where $p(\partial)$ is nonzero and assumed to be irreducible (see Section 4 below for more details). As $r(\partial)$ is in $\mathcal{P}$, it is an $\mathcal{A}$-linear combination of the $\ell$ rows $r_1(\partial), r_2(\partial), \ldots ,r_\ell(\partial)$ of $P(\partial)$, say $r(\partial) = b_1(\partial)r_1(\partial) + \ldots + b_\ell(\partial) r_\ell(\partial)$. Clearly the $b_j(\partial)$ are not all zero, and are also not all divisible by $p(\partial)$, because otherwise it would imply that $x(\partial)$ belongs to $\mathcal{P}$, contrary to its choice. Thus, without loss of generality, let $b_1(\partial)$ be nonzero and not divisible by $p(\partial)$. Now let $B(\partial)$ be the following $\ell \times \ell$ matrix: its first row is $(b_1(\partial), \ldots ,b_\ell(\partial))$, it has 1 in the $(j,j)$ entries, $2 \leqslant j \leqslant \ell$, and all other entries are 0. Then the product $B(\partial)P(\partial)$ is an $\ell \times k$ matrix whose first row is $r(\partial)$ and whose other rows are the rows $r_2(\partial), \ldots ,r_\ell(\partial)$ of $P(\partial)$. By construction of $B(\partial)$, its determinant $b_1(\partial)$ is not divisible by $p(\partial)$, and as every generator of the $\ell$-th determinantal ideal of $B(\partial)P(\partial)$ is divisible by $p(\partial)$, it follows that every generator of the $\ell$-th determinantal ideal $\frak{i}_\ell$ of $P(\partial)$ is also divisible by $p(\partial)$. This implies that $\rk(P(\lambda)) < \ell$ at points $\lambda$ in $\mathbb{C}^n$ where $p(\lambda) = 0$, an algebraic variety of dimension $n-1$.
\hspace*{\fill}$\square$\\

\noindent Remark 3.1: A remark similar to the one after Proposition 2.1 is relevant here (and elsewhere), namely that $\frak{i}_\ell$ is the $(k-\ell)$-th Fitting ideal of $\mathcal{A}^k/\mathcal{P}$. Thus, the assumption that $\frak{i}_\ell \neq 0$ implies that the $\ell$ rows of $P(\partial)$ is a minimum set of generators for $\mathcal{P}$ (in conformity with Definition 2.1).\\

\noindent Remark 3.2: A deceptively short proof of a related result on the controllability of discrete systems defined on $\mathbb{Z}^n$, using facts about Cohen-Macaulay rings, appears in Wood \cite{woo}, namely Corollaries 3.8 and 7.9. The longer proof above is elementary; besides the proof itself finds use in the next section.\\

\noindent Remark 3.3: Oberst shows in \cite{o}, that assuming $\frak{i}_\ell \neq 0$, the dimension of the variety $\mathcal{V}(\frak{i}_\ell)$ is less than or equal to $n-2$ if and only if the localization $\tor (\mathcal{A}^k/\mathcal{P})_\frak{p}$ at height 1 primes $\frak{p}$ is equal to 0, where $\tor(\mathcal{A}^k/\mathcal{P})$ denotes the submodule of torsion elements of $\mathcal{A}^k/\mathcal{P}$. The above theorem is a stronger statement, namely that assuming $\frak{i}_\ell \neq 0$, the dimension of $\mathcal{V}(\frak{i}_\ell)$ is less than or equal to $n-2$ if and only if $\tor (\mathcal{A}^k/\mathcal{P})$ is itself equal to 0. The apparent inconsistency  is resolved later in Remark 4.1 below. \\

The condition that $\frak{i}_\ell$ be nonzero in the above theorem, and remarks, admits an elementary interpretation.

\begin{proposition} Let $\mathcal{P}$ be a submodule of $\mathcal{A}^k$, and let $P(\partial)$ be any $\ell \times k$ matrix whose $\ell$ rows generate $\mathcal{P}$. Then $\mathcal{P}$ is free if its $\ell$-th determinantal ideal $\frak{i}_\ell$ is nonzero. Conversely, suppose that the $\ell$ rows of $P(\partial)$ is a minimum set of generators for $\mathcal{P}$. Then $\frak{i}_\ell$ is nonzero if $\mathcal{P}$ is free. Thus, the Hautus test of Theorem 3.1 is valid for behaviors defined by free submodules of $\mathcal{A}^k$.
\end{proposition}

\noindent Proof: Suppose $\mathcal{P}$ is not free; then there is a non-trivial relation $a_1r_{j_1} + \ldots + a_mr_{j_m} = 0$ between some of the $m$ ($ \leq \ell$) rows of $P(\partial)$, where by assumption $a_1, \dots , a_m$ are all nonzero. Let $P_1(\partial)$ be the $m \times k$ submatrix of $P(\partial)$ whose $m$ rows are $r_{j_1}, \ldots ,r_{j_m}$. Then the $m$-th determinantal ideal of $P_1(\partial)$, and hence also the $\ell$-th determinantal ideal of $P(\partial)$, equals 0.

Conversely, suppose that $\frak{i}_\ell = 0$. This means that after localizing at the 0 ideal, i.e over the function field $\mathcal{K} = \mathbb{C}(\partial_1, \ldots ,\partial_n)$, the $\ell$ rows of $P(\partial)$ span a subspace $\bar{\mathcal{P}}$ of $\mathcal{K}^k$, whose projections to $\mathcal{K}^\ell$, given by choosing $\ell$ of the $k$ coordinates, are all subspaces of dimension strictly less than $\ell$. If $k = \ell$, then the rows of $P(\partial)$ are $\mathcal{K}$-linearly dependent, hence there is an $\mathcal{A}$-relation between the rows of $P(\partial)$, and $\mathcal{P}$ is not free (as $\ell$ is the minimum number of elements needed to generate it). Otherwise, there are ${\tiny \left(\begin{array}{c}k \\\ell\end{array}\right)} > \ell$ such projections,  
hence $\bar{\mathcal{P}}$ must be of dimension strictly less than $\ell$, so that again, by definition of $\ell$, $\mathcal{P}$ is not free.
\hspace*{\fill}$\square$\\

The following example demonstrates that it is essential to assume that $\mathcal{P}$ be free in Theorem 3.1 (equivalently, that its $\ell$-th determinantal ideal be nonzero).\\

\noindent Example: Consider the submodules $\mathcal{P}_1$ and $\mathcal{P}_2$ of $\mathcal{A}^3$ generated by the matrices
\[
P_1(\partial) = \left(\begin{array}{ccc}0 & -\partial_3 & \partial_2 \\\partial_3 & 0 & -\partial_1 \\-\partial_2 & \partial_1 & 0\end{array}\right)
\hspace{3mm}  {\rm and} \hspace {4mm} P_2(\partial) = \left(\begin{array}{ccc}0 & -\partial_3 & \partial_2 \\\partial_3 & 0 & -\partial_1 \\-\partial_1 \partial_2 & \partial_1^2 & 0\end{array}\right)
\]
Here $k = \ell = 3$. The determinants of $P_1(\partial)$ and $P_2(\partial)$ are both 0, i.e. $\frak{i}_\ell = 0$, yet while $\mathcal{A}^3/\mathcal{P}_1$ is torsion free, $\mathcal{A}^3/\mathcal{P}_2$ has torsion. Hence, $\mathcal{B}_\mathcal{F}(\mathcal{P}_1)$ is controllable, whereas $\mathcal{B}_\mathcal{F}(\mathcal{P}_2)$ is not.
\hspace*{\fill}$\square$\\

The case when $\frak{i}_\ell$ equals 0, which is not answered by Theorem 3.1, is however not {\it generic} in a sense made precise later in the paper. The question now is whether there is a Hautus type test in this case. For instance, is controllability then determined by the first nonzero Fitting ideal of $\mathcal{A}^k/\mathcal{P}$?  

Consider the presentation
\[
\mathcal{A}^\ell \stackrel{P^t(\partial)}{\longrightarrow} \mathcal{A}^k \stackrel{\pi}{\longrightarrow} \mathcal{A}^k/\mathcal{P} \rightarrow 0
\] 
of $\mathcal{A}^k/\mathcal{P}$. If it is torsion free,  then it injects into a free module, say $\mathcal{A}^k/\mathcal{P} \stackrel{i}{\hookrightarrow}\mathcal{A}^{k_1}$, and in which case
\[
\mathcal{A}^\ell \stackrel{P^t(\partial)}{\longrightarrow} \mathcal{A}^k \stackrel{i \circ \pi}{\longrightarrow} \mathcal{A}^{k_1}
\]
is an exact sequence of free modules. This implies that $\dep(\mathcal{I}(P^t(\partial))) \geq 2$ (for instance, Theorem 20.9 of Eisenbud \cite{e}) where $\mathcal{I}(P^t(\partial))$ is the first nonzero Fitting ideal of $\mathcal{A}^k/\mathcal{P}$. As $\mathcal{A}$ is Cohen-Macaulay, depth equals codimension, and hence the dimension of the variety of $\mathcal{I}(P^t(\partial)) \leq n-2$. Thus a necessary condition for the controllability of $\mathcal{B}_\mathcal{F}(\mathcal{P})$ is that the dimension of the first nonzero Fitting ideal be less than or equal to $n-2$. This bound on the dimension cannot be a sufficient condition for controllability, because if $\mathcal{A}^k/\mathcal{P}$ is not torsion free, then the morphism $P^t(\partial)$ above fits into some exact sequence 
\[
\mathcal{A}^{\ell_1} \longrightarrow \mathcal{A}^\ell \stackrel{P^t(\partial)}{\longrightarrow} \mathcal{A}^k
\] 
and now, again by the theorem quoted above, $\dep(\mathcal{I}(P^t(\partial))) \geq 1$, so that a necessary condition for non-controllability of $\mathcal{B}_\mathcal{F}(\mathcal{P})$ is that the dimension of the first nonzero Fitting ideal be less than or equal to $n-1$. However this bound does not preclude it being less than or equal to $n-2$ as the following calculation indicates! \\

\noindent Example: The first nonzero Fitting ideals of $\mathcal{A}^3/\mathcal{P}_1$ and $\mathcal{A}^3/\mathcal{P}_2$ in the example above, are the first Fitting ideals, generated by the minors of the $2 \times 2$ submatrices of $P_1(\partial)$ and $P_2(\partial)$. These are $\mathcal{I}(P_1(\partial)) = (\partial_1^2, \partial_2^2, \partial_3^2, \partial_1\partial_2, \partial_1\partial_3, \partial_2\partial_3)$ and $\mathcal{I}(P_2(\partial)) = (\partial_1^3, \partial_3^2, \partial_1\partial_3, \partial_2\partial_3, \partial_1\partial_2^2, \partial_1^2\partial_2)$, whose radicals are $(\partial_1, \partial_2, \partial_3)$ and $(\partial_1, \partial_3)$ respectively. Both their dimensions are less than or equal to 1 (= 3-2), even though $\mathcal{A}^3/\mathcal{P}_2$ has torsion.
\hspace*{\fill}$\square$\\

Thus, there does not seem to be a straightforward generalization of the classical Hautus test when the ideal $\frak{i}_\ell$ equals 0, which is both necessary and sufficient. Theorem 3.1, however, does exhibit a `perfect' generalization of the classical Hautus test when this ideal is nonzero (and whose proof is also elementary). For another approach to this question, see Lomadze \cite{l1}.\\

The condition $\frak{i}_\ell \neq 0$ of Theorem 3.1 is however superfluous when $n = 1, 2$:\\

\noindent $n = 1$: Now every submodule $\mathcal{P}$ of $(\mathbb{C}[\frac{d}{dt}])^k$ is free (being torsion free). Moreover, algebraic subsets of dimension $n-2$ are empty. Hence, Theorem 3.1 specializes without qualification to the classical Hautus test, and to Proposition 2.1.\\

\noindent $n = 2$: Assume now that $\mathcal{A} = \mathbb{C}[\partial_1, \partial_2]$, and that the behavior $\mathcal{B}_\mathcal{F}(\mathcal{P})$ of a submodule $\mathcal{P}$ of $\mathcal{A}^k$ is controllable. Equivalently, $\mathcal{A}^k/\mathcal{P}$ is torsion free, and then it injects into a free module $\mathcal{A}^{k_1}$ as in the above presentation. This implies that in a minimal resolution
\[
\rightarrow \mathcal{A}^{\ell_1} \longrightarrow \mathcal{A}^\ell \stackrel{P^t(\partial)}{\longrightarrow} \mathcal{A}^k \longrightarrow \mathcal{A}^{k_1} \longrightarrow \mathcal{A}^{k_1}/(\mathcal{A}^k/\mathcal{P}) \rightarrow 0
\]
$\mathcal{A}^{\ell_1} = 0$, because the global dimension of the ring $\mathcal{A}$ equals 2. Thus the morphism $P^t(\partial)$ is injective, and this implies that its image $\mathcal{P}$ is free (see also Corollary 4 of \cite{ps}). Every controllable $2-D$ behavior is therefore given by a free submodule, so that the corresponding ideal $\frak{i}_\ell$ is always nonzero. Theorem 3.1 then implies the following result of Wood (\cite{woo}) and Zerz (\cite{z}) on $2-D$ behaviors:

\begin{corr} A $2-D$ behavior $\mathcal{B}_\mathcal{F}(\mathcal{P})$ is controllable if and only if $P(\lambda)$ drops rank at most at a finite number of points in $\mathbb{C}^2$ (where the $\ell$ rows of $P(\partial)$ is a minimum set of generators for $\mathcal{P}$).
\end{corr} 

\noindent Proof: By the above discussion it is unnecessary to assume that $\frak{i}_\ell \neq 0$, so that by Theorem 3.1, $\mathcal{B}_\mathcal{F}(\mathcal{P})$ is controllable if and only if the dimension of the variety of $\frak{i}_\ell$ is 0. The finite set of points of this variety is precisely where $P(\lambda)$ drops rank.
\hspace*{\fill}$\square$\\

The following example explains a classical terminology: \\

\noindent Example: Consider the special case of a {\it scalar} system given by $p_1(\partial)f = p_2(\partial)g$, i.e. the behavior defined by the kernel  of the map 
\[
\begin{array}{lccc}
\left(\begin{array}{cc} p_1(\partial), & -p_2(\partial)\end{array}\right): & \mathcal{F}^2 & \longrightarrow & \mathcal{F} \\ 
& (f,g) & \mapsto & (p_1(\partial),-p_2(\partial))\left(\begin{array}{c}f \\g\end{array}\right)
\end{array}
\]
The $\ell$-th determinantal ideal, here $\ell = 1$, is the ideal $\frak{i}_1 = (p_1,p_2)$ generated by $p_1$ and $p_2$. Theorem 3.1 asserts that this behavior is controllable if and only if there is no non-constant $p$ which divides both $p_1$ and $p_2$. This is the classical {\it pole-zero cancellation} criterion, a theorem in the behavioral setting (Willems \cite{w}). 
\hspace*{\fill}$\square$\\

This example suggests the following definition: 

\begin{defn} Let $\ell$ be the minimum number of elements needed to generate $\mathcal{P}$. Then the ideal $\frak{i}_\ell$ generated by the determinants of the $\ell \times \ell$ submatrices of (any) $P(\partial)$ whose $\ell$ rows generate $\mathcal{P}$, and its variety $\mathcal{V}(\frak{i}_\ell)$ (in $\mathbb{C}^n$), are called the {\it cancellation ideal} and the {\it cancellation variety} of the behavior $\mathcal{B}_\mathcal{F}(\mathcal{P})$ respectively. 
\end{defn}

Thus, Theorem 3.1 states that a behavior, whose cancellation variety is not all of $\mathbb{C}^n$, is controllable if and only if this variety is of dimension $\leqslant n-2$.

\begin{defn} Again let $\ell$ be the minimum number of elements needed to generate $\mathcal{P} \subset \mathcal{A}^k$. Its behavior $\mathcal{B}_\mathcal{F}(\mathcal{P})$ is said to be under-determined if $\ell \leqslant k$, and strictly under-determined if $\ell < k$. It is over-determined if $\ell \geqslant k$, and now it is the $k$-th determinantal ideal $\frak{i}_k$, generated by the determinants of $k \times k$ submatrices of (any) $P(\partial)$, whose $\ell$ rows generate $\mathcal{P}$, that is relevant: $\frak{i}_k$ is the 0-th Fitting ideal of $\mathcal{A}^k/\mathcal{P}$; it is the characteristic ideal of $\mathcal{B}_\mathcal{F}(\mathcal{P})$ (or of $\mathcal{P}$), and $\mathcal{V}(\frak{i}_k)$ is its characteristic variety \cite{o,sc}. Finally, a behavior is said to be square if $\ell = k$ (i.e. when it is both under-determined and over-determined).
\end{defn}

\noindent Remark 3.4: In the language that Oberst uses in \cite{o}, the points $\lambda$ in $\mathbb{C}^n$ where the rank of $P(\lambda)$ is strictly less than $\min\{\ell, k\}$ is the variety of {\it rank singularities}. It is more convenient for the purposes here to treat  the two cases separately, namely the cancellation variety when $\frak{i}_\ell\neq 0$, and the characteristic variety when $\frak{i}_k \neq 0$. \\

Theorem 3.1 implies that a behavior whose cancellation ideal $\frak{i}_\ell$ is nonzero, and which is controllable, must necessarily be strictly under-determined, the only exception being a square behavior whose cancellation ideal equals $\mathcal{A}$ (and then the behavior is 0). For if $\ell = k$, then $\frak{i}_\ell$ is a principal ideal not equal to 0 or $(1)$, hence $\mathcal{V}(\frak{i}_\ell)$ is of dimension $n-1$ and the behavior is not controllable, contradicting the assumption. The theorem below is the more general statement. It is included here even though it is well known, for instance \cite{o,ps,pq,woo,z}, in order to provide a convenient point of reference for results on genericity in Section 6. 

\begin{thm} Let $\mathcal{P}$ be a submodule of $\mathcal{A}^k$, and let 
$P(\partial)$ be any $\ell \times k$ matrix whose rows generate $\mathcal{P}$. Suppose that $\frak{i}_k$ is not the zero ideal (so that $\ell \geqslant k$). Then the behavior $\mathcal{B}_\mathcal{F}(\mathcal{P})$ is over-determined, and is not controllable unless $\frak{i}_k = \mathcal{A}$ (when the behavior equals 0).
\end{thm}

\noindent Proof: If $\mathcal{B}_\mathcal{F}(\mathcal{P})$ were not over-determined, then $\frak{i}_k$ would be the zero ideal. 

Now let $d$ be the determinant of a $k \times k$ submatrix $D(\partial)$ of $P(\partial)$, i.e. a generator of $\frak{i}_k$. Let $D'(\partial)$ be the matrix adjoint to $D(\partial)$, so that $D'D$ is the $k \times k$ diagonal matrix whose diagonal elements are all equal to $d$. But each row of $D'D$ is in $\mathcal{P}$, and thus it follows that if some element $a$ of $\mathcal{A}$ is in $\frak{i}_k$, then there is a $k \times k$ diagonal matrix whose diagonal entries are this $a$, and whose rows are in $\mathcal{P}$. Thus $\frak{i}_k = \mathcal{A}$ if and only if $\mathcal{P} = \mathcal{A}^k$. Its behavior then is 0 (and is trivially controllable).

Assume now that $\frak{i}_k \neq \mathcal{A}$, so that $\mathcal{P}$ is strictly contained in $\mathcal{A}^k$. As the rows of $D'D$ above are $d\cdot e_j, j = 1, \dots , k, ~e_j = (0, \dots , 1, \ldots , 0)$ -- 1 in the $j$-th place -- it also follows that $d$ is in the annihilator $\ann(\mathcal{A}^k/\mathcal{P})$ of $\mathcal{A}^k/\mathcal{P}$. This implies that 
\begin{equation}
\frak{i}_k \subseteq \ann(\mathcal{A}^k/\mathcal{P})
\end{equation}
and hence, as $\frak{i}_k$ is assumed to be nonzero, that the set of torsion elements of $\mathcal{A}^k/\mathcal{P}$ is nonzero. Then, by Theorem 1.2, it follows that $\mathcal{B}_\mathcal{F}(\mathcal{P})$ is not controllable.
\hspace*{\fill}$\square$\\

\noindent Remark 3.5: Consider a square behavior, i.e. $\ell = k$, so that  the assumption $\frak{i}_k \neq 0$ in Theorem 3.2 is equivalent to $\frak{i}_\ell \neq 0$ in Theorem 3.1. But then $\frak{i}_\ell$ is a (nonzero) principal ideal, its variety is of dimension $n-1$ (unless $\frak{i}_\ell = (1)$ when its variety is empty), and hence $\mathcal{B}_\mathcal{F}(\mathcal{P})$ is not controllable also by Theorem 3.1. Thus, Theorems 3.1 and 3.2 coincide when they are both the case. This again shows that a behavior whose cancellation ideal $\frak{i}_\ell$ is nonzero, and which is controllable, must necessarily be strictly under-determined (unless $P(\partial)$ is the identity matrix, i.e. $\mathcal{P} = \mathcal{A}^k$; then $\frak{i}_\ell = (1)$, and the behavior is 0).\\

In particular:

\begin{corr} Any controllable $2-D$ (nonzero) behavior is strictly under-determined.
\end{corr}

\noindent Proof: It is unnecessary to {\it assume} that $\frak{i}_\ell \neq 0$ when $n = 2$.
\hspace*{\fill}$\square$\\

\noindent Example: This statement is of course not true when $n \geq 3$. For instance, the behavior $\mathcal{B}_\mathcal{F}(\mathcal{P}_1)$ in the above examples is controllable ($P_1(\partial)$ is the $\curl$ operator), yet is not strictly under-determined.
\hspace*{\fill}$\square$

An over-determined system such that $\frak{i}_k \neq 0$ is also called an {\it autonomous} system; \cite{o,ps} describe its basic properties. Some of 
these properties follow directly from the Hautus test and are collected together in the last section to illustrate the test's importance, just as the classical Hautus test is important for lumped systems. This last section also studies the degenerate case when $\frak{i}_k = 0$. 

\vspace{.36cm}

\section{The structure of the cancellation variety} This section collects a few results about the cancellation ideal of a behavior, and its variety. They will be used in Section 6 on genericity, and also in Section 5 to obtain Hautus tests for controllability over spaces other than $\mathcal{D}'$ and $\mathcal{C}^\infty$. These results follow from the {\it proof} of Theorem 3.1.

Recollect that a prime ideal $\frak{p}$ is an associated prime of an $\mathcal{A}$-module $\mathcal{M}$ if it is equal to the annihilator $\ann(x)$ of some element $x$ of $\mathcal{M}$. An element $a$ of $\mathcal{A}$ is a zero divisor for $\mathcal{M}$ if there is a nonzero $x$ in $\mathcal{M}$ such that $ax = 0$. The maximal elements of the family $\{\ann(x) | 0 \neq x \in \mathcal{M}\}$ are associated primes of $\mathcal{M}$, hence the union of all the associated primes is the set of all the zero divisors for $\mathcal{M}$ (for instance \cite{e}). It then follows that $\mathcal{M}$ is torsion free if and only if 0 is its only associated prime, and hence $\mathcal{B}_\mathcal{F}(\mathcal{P})$ is controllable if and only if 0 is the only associated prime of $\mathcal{A}^k/\mathcal{P}$.

Returning to the case at hand, namely when the module $\mathcal{M}$ equals $\mathcal{A}^k/\mathcal{P}$ with $\mathcal{P}$ a free submodule of $\mathcal{A}^k$ (equivalently, whose cancellation ideal $\frak{i}_\ell$ is nonzero), the first half of Theorem 3.1 demonstrates that an irreducible $p(\partial)$ which divides every generator of $\frak{i}_\ell$ is a zero divisor for $\mathcal{A}^k/\mathcal{P}$. 

Conversely, suppose that $a(\partial)$ is a zero divisor for $\mathcal{A}^k/\mathcal{P}$, and that $x$ is an element of $\mathcal{A}^k \setminus \mathcal{P}$ such that $ax \in \mathcal{P}$. Because $\mathcal{A}$ here is a UFD, $a$ is a unique product of irreducible factors, and it follows (as observed in the second half of the proof of Theorem 3.1) that there is an $x'$ in $\mathcal{A}^k \setminus \mathcal{P}$ such that $px' \in \mathcal{P}$ for some irreducible factor $p(\partial)$ of $a(\partial)$. In other words, the maximal elements of the family of principal ideals generated by zero divisors of $\mathcal{A}^k/\mathcal{P}$ are the principal ideals generated by irreducible zero divisors. Furthermore, by the second half of the proof of Theorem 3.1, every irreducible zero divisor divides every generator of $\frak{i}_\ell$.  Together with the above paragraph, this shows that the set of irreducible zero divisors for $\mathcal{A}^k/\mathcal{P}$ is precisely the set of irreducible common factors of the ${\tiny \left(\begin{array}{c}k \\\ell\end{array}\right)}$ generators of the cancellation ideal $\frak{i}_\ell$. As the latter set is finite, it follows that the set of irreducible zero divisors for $\mathcal{A}^k/\mathcal{P}$ is a finite set.

This in turn implies that the nonzero associated primes of $\mathcal{A}^k/\mathcal{P}$ are precisely the principal ideals generated by the irreducible zero divisors described above. For suppose that some associated prime $\frak{p}$ is not principal. Let $p_1$ be an irreducible element of $\frak{p}$ 
such that its degree is minimum amongst all elements in $\frak{p}$. Let $p_2$ be an irreducible element in $\frak{p} \setminus (p_1)$, again of minimum degree. Then it is elementary that $p_1 + \alpha p_2$ is irreducible for all 
$\alpha$ in $\mathbb{C}$. These infinite number of ireducible elements are all in $\frak{p}$, and hence are all zero divisors for $\mathcal{A}^k/\mathcal{P}$, contradicting the assertion above. This establishes the following theorem:

\begin{thm} Let $\mathcal{P}$ be a free submodule of $\mathcal{A}^k$ (that is, whose cancellation ideal $\frak{i}_\ell$ is nonzero). Then the nonzero associated primes of $\mathcal{A}^k/\mathcal{P}$ are the principal ideals generated by the irreducible $p(\partial)$ that divide every generator of $\frak{i}_\ell$.
\hspace*{\fill}$\square$
\end{thm}

\noindent Aside: A part of the above theorem, that the associated primes are principal, also follows from general considerations. Let $\frak{m}$ be a maximal ideal containing an associated prime $\frak{p}$ of $\mathcal{A}^k/\mathcal{P}$. Localizing at $\frak{m}$, $\mathcal{P}_\frak{m}$ is a free submodule of $\mathcal{A}^k_\frak{m}$, hence the following sequence
\[
0 \rightarrow \mathcal{A}^\ell_\frak{m} \longrightarrow \mathcal{A}^k_\frak{m} \longrightarrow \mathcal{A}^k_\frak{m}/\mathcal{P}_\frak{m} \rightarrow 0
\]
is exact. This implies that the projective dimension $\pd (\mathcal{A}^k_\frak{m}/\mathcal{P}_\frak{m})$ is equal to 1. By Auslander-Buchsbaum (for instance, Theorem 19.9 in \cite{e}), $\dep (\mathcal{A}^k_\frak{m}/\mathcal{P}_\frak{m}) = \dep (\mathcal{A}_\frak{m}) - \pd (\mathcal{A}^k_\frak{m}/\mathcal{P}_\frak{m})$, which equals $n-1$ (as $\mathcal{A}_\frak{m}$ is Cohen-Macaulay, its depth is equal to $n$, its dimension). On the other hand, $\dep (\mathcal{A}^k_\frak{m}/\mathcal{P}_\frak{m}) \leq \dim (\mathcal{A}_\frak{m}/\frak{p})$ (Proposition 18.2 in \cite{e} for instance), hence it follows that $\height (\frak{p}) = \dim (\mathcal{A}_\frak{m}) - \dim (\mathcal{A}_\frak{m}/\frak{p}) \leq 1$, so that $\frak{p}$ is a principal ideal.
\hspace*{\fill}$\square$\\

The following corollary is now immediate.

\begin{corr} The affine varieties of the nonzero associated primes of $\mathcal{A}^k/\mathcal{P}$, when the cancellation ideal $\frak{i}_\ell$ of $\mathcal{P}$ is nonzero, are the irreducible components of dimension $n-1$ of the cancellation variety $\mathcal{V}(\frak{i}_\ell)$. 
\hspace*{\fill}$\square$
\end{corr}

\noindent A remark on the controllable-uncontrollable decomposition of a behavior (see also \cite{o}): Suppose $\mathcal{B}_\mathcal{F}(\mathcal{P})$ is a general behavior determined by a (not necessarily free) submodule $\mathcal{P}$ of $\mathcal{A}^k$. Let $\mathcal{P}'$ be the submodule $\{x \in \mathcal{A}^k | ax \in \mathcal{P}, a \neq 0\}$. Then $\mathcal{P}'$ contains $\mathcal{P}$, and the quotient $\mathcal{P}'/\mathcal{P}$ is the submodule $\tor (\mathcal{A}^k/\mathcal{P})$ of torsion elements of $\mathcal{A}^k/\mathcal{P}$. The following sequence 
\[
0 \rightarrow \mathcal{P}'/\mathcal{P} \longrightarrow \mathcal{A}^k/\mathcal{P} \longrightarrow \mathcal{A}^k/\mathcal{P}' \rightarrow 0
\]
is exact, where $\mathcal{A}^k/\mathcal{P}'$ is torsion free. In general, given a short exact sequence of $\mathcal{A}$-modules, the associated primes of the middle term is contained in the union of the associated primes of the other two modules \cite{e}. Here however, it is clear, that there is equality.

\begin{lemma} $\ass (\mathcal{A}^k/\mathcal{P}) = \ass (\mathcal{P}'/\mathcal{P}) \bigcup \ass (\mathcal{A}^k/\mathcal{P}')$. Hence, if $\mathcal{P} \subsetneq\mathcal{P}' \subsetneq \mathcal{A}^k$, then $\ass (\mathcal{A}^k/\mathcal{P}') = \{0\}$ and $\ass (\mathcal{P}'/\mathcal{P})$ is the set of all the nonzero associated primes of $\mathcal{A}^k/\mathcal{P}$.
\hspace*{\fill}$\square$
\end{lemma}

\noindent Remark 4.1: Under the assumption $\frak{i}_\ell \neq 0$, the nonzero associated primes of $\mathcal{A}^k/\mathcal{P}$ are all of height 1 (Theorem 4.1 above). Thus the associated primes of $\tor(\mathcal{A}^k/\mathcal{P})$ are precisely these height 1 primes (the above lemma). Therefore, to say that the localization $\tor(\mathcal{A}^k/\mathcal{P})_\frak{p}$ equals 0 at height 1 primes is to say that the set of associated primes of $\tor(\mathcal{A}^k/\mathcal{P})$ is empty, which in turn is to say that $\tor(\mathcal{A}^k/\mathcal{P}) = 0$. This resolves the `inconsistency' discussed in Remark 3.3. \\

Applying the functor $\homo_\mathcal{A}(-, \mathcal{F})$, $\mathcal{F}$ either $\mathcal{C}^\infty$ or $\mathcal{D}'$, to the above sequence gives 
\begin{center}
$ 0 \rightarrow \mathcal{B}_\mathcal{F}(\mathcal{P}') \longrightarrow \mathcal{B}_\mathcal{F}(\mathcal{P}) \longrightarrow \homo_\mathcal{A}(\mathcal{P}'/\mathcal{P},~\mathcal{F}) \rightarrow 0 $
\end{center}
which is exact as $\mathcal{F}$ is an injective $\mathcal{A}$-module. As $\mathcal{A}^k/\mathcal{P}'$ is torsion free, the behavior $\mathcal{B}_\mathcal{F}(\mathcal{P}')$ is controllable, and is a sub-behavior of $\mathcal{B}_\mathcal{F}(\mathcal{P})$; in fact, it is the largest controllable sub-behavior of $\mathcal{B}_\mathcal{F}(\mathcal{P})$ in the sense that any other controllable sub-behavior is contained in it. The {\it quotient} behavior $ \homo_\mathcal{A}(\mathcal{P}'/\mathcal{P}, \mathcal{F}) \simeq \mathcal{B}_\mathcal{F}(\mathcal{P})/\mathcal{B}_\mathcal{F}(\mathcal{P}')$ is not only not controllable, it also does not contain any nonzero controllable sub-behaviors, and is {\it uncontrollable} by definition. Note that this uncontrollable behavior is not a sub-behavior of $\mathcal{B}_\mathcal{F}(\mathcal{P})$ unless the above sequences split. In the case of lumped behaviors, i.e. when $n = 1$, $\mathcal{A}^k/\mathcal{P}'$ is free being torsion free, hence the sequences split, and then this uncontrollable behavior can also be considered a sub-behavior of $\mathcal{B}_\mathcal{F}(\mathcal{P})$, \cite{w}. 

Thus, if $\mathcal{P} = \mathcal{P}'$, then $\mathcal{A}^k/\mathcal{P}$ is torsion free (equivalently, 0 is its only associated prime) and $\mathcal{B}_\mathcal{F}(\mathcal{P})$ is controllable, while if $\mathcal{P}' = \mathcal{A}^k$, then $\mathcal{A}^k/\mathcal{P}$ is a torsion module (equivalently, 0 is not an associated prime), and $\mathcal{B}_\mathcal{F}(\mathcal{P})$ does not contain any nonzero controllable behavior, hence is uncontrollable.

Recollect that the minimal elements of the set of associated primes of $\mathcal{M}$ coincide with the minimal elements of the support $\supp(\mathcal{M})$ of $\mathcal{M}$, where $\supp(\mathcal{M}) = \{\frak{p} ~| ~\mathcal{M}_\frak{p} \neq 0\}$. When $\mathcal{M}$ is finitely generated, these are also the minimal primes containing the annihilator $\ann(\mathcal{M})$ of $\mathcal{M}$. These observations immediately lead to a strengthening of Theorem 3.2.

\begin{thm} Let $\mathcal{P}$ be strictly contained in $\mathcal{A}^k$, and let $P(\partial)$ be any $\ell \times k$ matrix whose rows generate $\mathcal{P}$. Suppose that $\frak{i}_k$ is not the zero ideal (so that $\ell \geqslant k$). Then the over-determined system $\mathcal{B}_\mathcal{F}(\mathcal{P})$ is uncontrollable, i.e. it does not contain any nonzero controllable sub-behavior.
\end{thm}

\noindent Proof: Inclusion (3) in the proof of Theorem 3.2 implies as before that $\ann(\mathcal{A}^k/\mathcal{P})$ is nonzero, and hence that 0 is not an associated prime of $\mathcal{A}^k/\mathcal{P}$ (see also \cite{ps}).
\hspace*{\fill}$\square$\\

\vspace{.3cm}

\section {The Hautus test for other signal spaces} The above description of the associated primes of $\mathcal{A}^k/\mathcal{P}$ when $\frak{i}_\ell$ is nonzero (namely Theorem 4.1 and Corollary 4.1), leads to a Hautus test for other signal spaces $\mathcal{F}$. It will be convenient in this section to use the notation $D_j = \frac{1}{\imath}\frac{\partial}{\partial_j}$, so that the ring $\mathcal{A}$ of partial differential operators is also the polynomial ring $\mathbb{C}[D_1, \dots ,D_n]$. 

The first step in obtaining a Hautus test for other spaces $\mathcal{F}$ is to reformulate the `patching' definition of controllability of Definition 1.1 so as to be meaningful for these $\mathcal{F}$. This is easily done using the `cutoff' formulation of controllability in \cite{ps}. 

\begin{lemma} \cite{ps} A behavior $\mathcal{B}$ in $\mathcal{D}'$ or $\mathcal{C}^\infty$ is controllable if and only if for any subset $U$ of $\mathbb{R}^n$ and any open subset $O$ that contains the closure of $U$, there is for every $f$ in $\mathcal{B}$ an element $f_c$ in $\mathcal{B}$ that equals $f$ on some neighborhood of $U$ and equals 0 in the complement of $O$. Such an $f_c$ is said to be a cutoff of $f$ with respect to $U$ and $O$.
\hspace*{\fill}$\square$
\end{lemma}

This cutoff formulation of controllability cannot be carried over {\it in toto} to an arbitrary $\mathcal{F}$ because it might be that with respect to some $U$ and $O$, there is no cutoff $f_c$ of an $f$ in $\mathcal{B}$, even in $\mathcal{F}^k$. \\

\noindent Example: Let $\mathcal{F}$ be the space $\mathcal{S}'$ of temperate distributions on $\mathbb{R}^2$. Let $U = \{(x,y) \in \mathbb{R}^2 | y < 0\}$ and $O = \{(x,y) \in \mathbb{R}^2 | y < e^{-x^2}\}$. Then no cutoff of a constant function with respect to this $U$ and $O$ can be in $\mathcal{S}'$ (the derivatives of any such cutoff would not be tempered). 
\hspace*{\fill}$\square$\\

Thus the only possible way to generalize the cutoff formulation of controllability to an arbitrary $\mathcal{F}$ is the following:

\begin{defn} Let the signal space be an $\mathcal{A}$-submodule $\mathcal{F}$ of $\mathcal{D}'$. A subset $U$ of $\mathbb{R}^n$ and an open set $O$ that contains its closure are said to be admissible with respect to $\mathcal{F}$ if every $f$ in it admits a cutoff $f_c$ in $\mathcal{F}$ with respect to $U$ and $O$.

A behavior $\mathcal{B}$ in $\mathcal{F}$ is said to be controllable if for any $U$ and $O$ admissible with respect to $\mathcal{F}$, every $f$ in $\mathcal{B}$ admits a cutoff $f_c$ in $\mathcal{B}$ with respect to $U$ and $O$.
\end{defn}
    
Thus the subsets $U$ and $O$ of $\mathbb{R}^2$ in the above example are {\it not} admissible with respect to $\mathcal{S}'$. Some examples of admissible $U$ and $O$ are given in the following proposition.

\begin{proposition} (i) Let $\mathcal{F}$ be $\mathcal{D}'$ or $\mathcal{C}^\infty$. Then any pair $U$ and $O$ as in the above definition is admissible (for which reason, no qualification was ever necessary thus far in the paper).

\vspace{1mm}
\noindent (ii) Suppose that the distance between $U$ and the complement $O^c$ of $O$ is bounded away from 0, i.e. $\|x-y\| \geqslant \epsilon > 0$ for all $x \in U, y \in O^c$ (for instance, suppose the boundary of $U$ is compact). Then such a pair is admissible with respect to the space $\mathcal{S}'$ of temperate distributions.
\end{proposition}

\noindent Proof: (i) is elementary, for let $\rho$ be any smooth function which equals 1 on some neighborhood $U'$ of $U$ in $O$, and 0 on the complement of $O$. Then for an element $f$ in $\mathcal{D}'$ or $\mathcal{C}^\infty$, $\rho f$ is the required cutoff.

\vspace{1mm}
(ii) Let $U'$ be an open subset containing $U$ and contained in $O$, such that the distance between $U'$ and the complement of $O$ is $\frac{\epsilon}{2}$. Let $\chi$ be the characteristic function of $U'$, and let $\kappa$ be a smooth `bump' function supported in the ball of radius $\frac{\epsilon}{4}$ centered at the origin, i.e. $\kappa$ is identically 1 in some smaller ball about the origin. Then the convolution $\rho = \kappa \star \chi$ is a smooth function that is identically 1 on $U'$, 0 on the complement of $O$, and such that all its derivatives are bounded on $\mathbb{R}^n$. Then for any $f$ in $\mathcal{S}'$, $\rho f$ is also in $\mathcal{S}'$.
\hspace*{\fill}$\square$\\

Consider now the space $\mathcal{S}'$ of temperate distributions on $\mathbb{R}^n$. It is an $\mathcal{A}$-submodule of $\mathcal{D}'$ which is injective, but not a cogenerator (for instance \cite{sg}). Necessary and sufficient conditions for controllability of a behavior in $\mathcal{S}'$ is the following result from \cite{sg}.

\begin{thm} Let $\mathcal{A} = \mathbb{C}[D_1, \ldots ,D_n]$, and let the signal space $\mathcal{F}$ be the space $\mathcal{S}'$ of temperate distributions on $\mathbb{R}^n$. Let $\mathcal{P}$ be a submodule of $\mathcal{A}^k$. Then the behavior $\mathcal{B}_{\mathcal{S}'}(\mathcal{P})$ is controllable (in the sense of Definition 6.1) if and only if the varieties of the nonzero associated primes of $\mathcal{A}^k/\mathcal{P}$ (in $\mathbb{C}^n$) do not contain real points.
\end{thm}

\noindent Proof: This statement is the translation of Theorem 3.1 (iii) of \cite{sg} when the ring of differential operators is written as polynomials in $D_j = \frac{1}{\imath}\partial_j$ (instead of in $\partial_j$), and gives necessary and sufficient conditions for a behavior, given as a kernel, to admit an image representation. Images are clearly controllable in the sense of Definition 5.1, for the image of an arbitrary cutoff, with respect to any admissible $U$ and $O$, is a cutoff of an element in the behavior that also  belongs to the behavior.    
\hspace*{\fill}$\square$\\

Combining this with Corollary 4.1 above, gives the analogue of Theorem 3.1 for behaviors in the space $\mathcal{S}'$.

\begin{thm} Let $\mathcal{P}$ be a free submodule of $\mathcal{A}^k$, and let $P(D), ~D=\frac{1}{\imath}\partial,$ be any $\ell \times k$ matrix whose $\ell$ rows is a basis for $\mathcal{P}$ (so that the $\ell$-th determinantal ideal $\frak{i}_\ell$ of $P(D)$ is not zero). Then $\mathcal{B}_{\mathcal{S}'}(\mathcal{P})$ is controllable if and only if $P(\lambda)$ has full row rank for all $\lambda$ in the complement of an algebraic variety in $\mathbb{C}^n$ of dimension $\leqslant n-1$, where every irreducible component of this variety, of dimension $n-1$, does not intersect $\mathbb{R}^n$. 
\end{thm}

\noindent Proof: By Corollary 4.1, if there is an irreducible variety of dimension $n-1$ with real points on which $P(\lambda)$ drops rank, then the principal ideal generated by the corresponding irreducible polynomial is an associated prime of $\mathcal{A}^k/\mathcal{P}$, whose variety therefore contains real points. By Theorem 5.1, $\mathcal{B}_{\mathcal{S}'}(\mathcal{P})$ is then not be controllable.  
\hspace*{\fill}$\square$\\

Similar considerations apply to the following spaces defined by periodic functions \cite{nps}. Let $\t$ be the torus $\mathbb{R}^n/2\pi \mathbb{Z}^n$. Let $\mathcal{C}^\infty(\t)$ be the space of smooth functions on $\t$, in other words, smooth functions on $\mathbb{R}^n$ periodic with respect to the lattice $2\pi \mathbb{Z}^n$. Let $\mathcal{C}^\infty(\t)_{\rm fin}$ be the $\mathcal{A}$-submodule of $\mathcal{C}^\infty(\t)$ consisting of those periodic functions whose Fourier series expansion is a finite sum. Let $\mathcal{F} = \mathcal{C}^\infty(\t)_{\rm fin}[x_1, \ldots, x_n]$ be the subalgebra of $\mathcal{C}^\infty(\mathbb{R}^n)$ obtained by adjoining the coordinate functions $x_1, \ldots, x_n$ to $\mathcal{C}^\infty(\t)_{\rm fin}$.  

For positive integers $N_1$ dividing $N_2$, the natural $\mathcal{A}$-module morphism $\mathcal{C}^\infty (\mathbb{R}^n/2\pi N_1\mathbb{Z}^n) \rightarrow  
\mathcal{C}^\infty (\mathbb{R}^n/2\pi N_2\mathbb{Z}^n)$ identifies the space of functions periodic with respect to the lattice $2\pi N_1\mathbb{Z}^n$ with a closed subspace of the space of functions periodic with respect to the larger lattice $2\pi N_2\mathbb{Z}^n$. Let $\mathcal{C}^\infty(\pt)$ be the direct limit $\underset{\rightarrow }{\lim} \ \mathcal{C}^\infty (\mathbb{R}^n/2\pi N\mathbb{Z}^n)$. (These functions live on a {\it protorus}, which is by definition the inverse limit of the tori $\mathbb{R}^n/2\pi N\mathbb{Z}^n$, see \cite{nps}.) As before let $\mathcal{C}^\infty(\pt)_{\rm fin}$ be those functions in $\mathcal{C}^\infty (\pt)$ whose Fourier expansion is a finite sum. Finally, let $\mathcal{C}^\infty(\pt)_{\rm fin}[x_1, \ldots ,x_n]$ be the subalgebra of $\mathcal{C}^\infty(\mathbb{R}^n)$ obtained by adjoining $x_1, \ldots ,x_n$ to $\mathcal{C}^\infty (\pt)_{\rm fin}$.

Necessary and sufficient conditions for a behavior, given as a kernel and to admit an image representation in these spaces, is the following theorem.

\begin{thm} (Theorem 4.2 in \cite{nps}) (i) Let $\mathcal{A} = \mathbb{C}[D_1, \ldots ,D_n]$, and let $\mathcal{F}=\mathcal{C}^\infty(\pt)_{\rm fin}[x_1,\dots ,x_n]$. Then the behavior $\mathcal{B}_\mathcal{F}(\mathcal{P})$ of a submodule $\mathcal{P} \subset \mathcal{A}^k$ is the image of some morphism, and hence controllable in the sense of Definition 5.1, if and only if the varieties of the nonzero associated primes of $\mathcal{A}^k/\mathcal{P}$ do not contain rational points.

(ii) Now let $\mathcal{F}=\mathcal{C}^\infty(\t)_{\rm fin}[x_1,\dots ,x_n]$. Then $\mathcal{B}_\mathcal{F}(\mathcal{P})$ is the image of some morphism, and hence controllable in the sense of Definition 6.1, if and only if the varieties of the nonzero associated primes of $\mathcal{A}^k/\mathcal{P}$ do not contain integral points.
\hspace*{\fill}$\square$
\end{thm} 

Corollary 4.1 now gives the analogues of Theorem 3.1.

\begin{thm} Let $\mathcal{P}$ be a free submodule of $\mathcal{A}^k$, and let $P(D), ~D=\frac{1}{\imath}\partial,$ be any $\ell \times k$ matrix whose $\ell$ rows is a basis for $\mathcal{P}$ (so that the $\ell$-th determinantal ideal $\frak{i}_\ell$ of $P(D)$ is not zero). 

Suppose $\mathcal{F} = \mathcal{C}^\infty(\pt)_{\rm fin}[x_1,\dots ,x_n]$ (respectively, $\mathcal{F} = \mathcal{C}^\infty(\t)_{\rm fin}[x_1,\dots ,x_n]$). Then $\mathcal{B}_\mathcal{F}(\mathcal{P})$ admits an image representation, and is therefore controllable, if and only if $P(\lambda)$ has full row rank for all $\lambda$ in the complement of an algebraic variety in $\mathbb{C}^n$ of dimension $\leqslant n-1$, where every irreducible component of this variety, of dimension $n-1$, does not intersect $\mathbb{Q}^n$ (respectively, does not intersect $\mathbb{Z}^n$).
\hspace*{\fill}$\square$
\end{thm}

\section{Genericity} The purpose of this section is to show that a generic strictly under-determined behavior is controllable, whereas a generic over-determined behavior is uncontrollable. 

These genericity results are with respect to the choice of a topology on the set of behaviors. The topology discussed below is dictated by considerations of {\em structured perturbations}, which is natural in this context. For instance, perturbations of an under-determined behavior should result only in under-determined behaviors (and similarly for over-determined behaviors). This is because we expect that perturbing a behavior, which is equivalent to perturbing the differential equations defining the behavior, should involve perturbing the coefficients that appear in the differential equations, and should {\em not} involve any change in the basic structure of the defining differential equations, such as a change in the number of equations, for instance. Thus, if the behavior is defined by a submodule $\mathcal{P} \subseteq \mathcal{A}^k$, i.e. if the behavior is $\mathcal{B}_\mathcal{F}(\mathcal{P})$, and if $\mathcal{P}$ can be generated by $\ell$ elements, then any perturbation of this behavior should be of the form $\mathcal{B}_\mathcal{F}(\mathcal{P}')$, where the perturbed submodule $\mathcal{P}'$ can also be generated by $\ell$ elements. This is a reasonable requirement because the equations defining the behavior would arise from physical laws governing the behavior, such as the Newton, Maxwell, or thermodynamic equations, and from structural configurations involving interconnections of sub-systems etc. (see Willems \cite{w}). Thus, these equations, when perturbed, would result in similar equations but with different coefficients.

This section is in two parts. In the first, genericity results are established for matrices $\mathcal{M}_{\ell,k}$ of size $\ell \times k$, with entries in $\mathcal{A} = \mathbb{C}[\partial_1, \ldots ,\partial_n]$. The above discussion naturally leads to two distinct cases, the strictly under-determined when $\ell < k$, and the over-determined when $\ell \geqslant k$. In the second part, these genericity results for matrices are used to obtain corresponding genericity results for behaviors. Here some care must be exercised, because in general, the rows of infinitely many matrices in $\mathcal{M}_{\ell,k}$ generate the same submodule of $\mathcal{A}^k$. The difficulty in this descent from matrices to behaviors is compounded by the fact that matrices with different number of rows, i.e. different $\ell$, can generate the same submodule. This difficulty is circumvented by working with the subset $\mathcal{R}_{\ell,k}$ of those matrices whose rows are a minimum set of generators for the submodules they generate (namely Definition 2.1).  It turns out that the complement of $\mathcal{R}_{\ell,k}$ is contained in a proper Zariski closed subset of $\mathcal{M}_{\ell,k}$, and is therefore a vanishingly thin subset. Restricting to $\mathcal{R}_{\ell,k}$ permits the definition of a topology on the set of behaviors in a consistent manner satisfying the requirements of structured perturbations explained above. Genericity results for behaviors with respect to this topology then follow from genericity results for $\mathcal{M}_{\ell,k}$. \\

Consider first the case of under-determined matrices $\mathcal{M}_{\ell,k}$ with $\ell \leqslant k$. Let $\mathcal{M}_{\ell, k}(d)$ be the subset of those matrices in $\mathcal{M}_{\ell,k}$ whose entries are all bounded in degree by $d$. For $d_1 < d_2$, there is a natural inclusion $\mathcal{M}_{\ell,k}(d_1) \hookrightarrow \mathcal{M}_{\ell,k}(d_2)$, and as $d$ tends to infinity, the direct limit of $\mathcal{M}_{\ell, k}(d)$ is $\mathcal{M}_{\ell,k}$. The lemma below states that the set of $\ell \times k$ matrices such that $\frak{i}_\ell \neq 0$ is Zariski open in $\mathcal{M}_{\ell, k}(d)$ for every $d$, and hence in the direct limit $\mathcal{M}_{\ell,k}$. Finally allowing $\ell$ to vary between 1 and $k$ shows that this condition on $\frak{i}_\ell$ is generically satisfied on the disjoint union $\mathcal{M}_{\ell \leq k} = \bigcup_{\ell = 1}^k \mathcal{M}_{\ell, k}$.

\begin{lemma} In the notation of the above paragraph, $\frak{i}_\ell \neq 0$ is satisfied by elements in a Zariski open subset of  $\mathcal{M}_{\ell \leq k}$.
\end{lemma}

\noindent Proof: Consider the set $\mathcal{M}_{\ell, k}(d)$, and let $P(\partial)$ be an element of it. Each entry of $P(\partial)$ is a polynomial in $n$ indeterminates, of degree at most d. There are ${\tiny \left(\begin{array}{c}n+d \\n\end{array}\right)}$ such monomials, and there are $\ell k$ entries, hence $P(\partial)$ is a point in $\mathbb{C}^{N(d)}$ where $N(d) = \ell k{{\tiny \left(\begin{array}{c}n+d \\n\end{array}\right)}}$ (and thus $N(d)$ is a polynomial in $d$ of degree $n$).

The determinant of an $\ell \times \ell$ submatrix of $P(\partial)$ in $\mathcal{M}_{\ell, k}(d)$ is bounded in degree by $\ell d$, and there are 
${\tiny \left(\begin{array}{c}n+\ell d \\n\end{array}\right)}$ such monomials. There are ${\tiny \left(\begin{array}{c}k \\\ell\end{array}\right)}$ many minors, hence the $\ell$-th determinantal ideal of $P(\partial)$ is a point in $\mathbb{C}^{M(d)}$, where $M(d) = {\tiny \left(\begin{array}{c}k \\\ell\end{array}\right)}{\tiny \left(\begin{array}{c}n+\ell d \\n\end{array}\right)}$ is again a polynomial in $d$ of degree $n$. Thus the $\ell$-th determinantal ideal is a morphism
\[
\frak{i}_\ell(d): \mathbb{C}^{N(d)} \longrightarrow \mathbb{C}^{M(d)}
\]
where $\frak{i}_\ell(d)$ coincides with $\frak{i}_\ell(d-1)$ on $\mathbb{C}^{N(d-1)} \hookrightarrow \mathbb{C}^{N(d)}$.
Each of the $M(d)$ components of the morphism $\frak{i}_\ell(d)$ is polynomial in the entries of the argument, and is therefore algebraic. This implies that $\frak{i}_\ell(d)^{-1}(0)$ is a proper Zariski closed subset, and hence that the set of $P(\partial)$ in $\mathcal{M}_{\ell, k}(d)$, for which $\frak{i}_\ell (d) \neq 0$, is Zariski open. 

For $d_1 < d_2$, $\mathcal{M}_{\ell, k}(d_1)$ injects into $\mathcal{M}_{\ell, k}(d_2)$, and the topology of the former is the topology it inherits as a subspace of the latter (and thus $\mathcal{M}_{\ell, k}$ is a strict direct limit of the directed system  $\{\mathcal{M}_{\ell, k}(d)\}_{d=1,2,...}$). The lemma now follows from the definition of the topology on the direct limit.  
\hspace*{\fill}$\square$\\

\begin{corr}The Hautus tests of Theorems 3.1, 5.2 and 5.4 are valid for a generic under-determined matrix, namely those belonging to a Zariski open set given by the condition $\frak{i}_\ell \neq 0$.
\hspace*{\fill}$\square$
\end{corr}

An argument similar to the proof of Lemma 6.1 gives:

\begin{lemma}Consider the set of strictly under-determined matrices $\mathcal{M}_{\ell < k} = \bigcup_{\ell = 1}^{k-1}\mathcal{M}_{\ell, k}$. Then the condition that two or more generators of $\frak{i}_\ell$ are nonzero is satisfied by elements belonging to a Zariski open subset.
\end{lemma}

\noindent Proof: The set of all elements in $\mathcal{M}_{\ell,k}(d)$ such that exactly one of the ${\tiny \left(\begin{array}{c}k \\\ell\end{array}\right)}$ minors is nonzero, is a closed subset of the open set $\mathbb{C}^{N(d)} \setminus \frak{i}_\ell (d)^{-1}(0)$. 
\hspace*{\fill}$\square$\\

Now let $p(\partial)$ be a nonconstant operator in $\mathcal{A}$ of degree $\delta$. Multiplication in $\mathcal{A}$ by $p(\partial)$ increases degree by $\delta$, and hence there is an injective morphism
\[
p(\partial): \mathbb{C}^{{\tiny \left(\begin{array}{c}n+d \\n\end{array}\right)}} \longrightarrow \mathbb{C}^{{\tiny \left(\begin{array}{c}n+d+\delta \\n\end{array}\right)}}
\]
where its domain is the space of all operators in $\mathcal{A}$ of degree less than or equal to $d$, and its codomain, the space of all operators of degree less than or equal to $d+\delta$.

\begin{lemma}The image of the above morphism $p(\partial)$ is contained in a proper Zariski closed subset of the codomain.
\end{lemma}

\noindent Proof: Let X be the closure of the image of the morphism $p(\partial)$. Then X is irreducible, and the morphism $p(\partial): \mathbb{C}^{{\tiny \left(\begin{array}{c}n+d \\n\end{array}\right)}} \rightarrow X$ is dominant. This implies that $\dim(X) = {\small \left(\begin{array}{c}n+d \\n\end{array}\right) < \left(\begin{array}{c}n+d+\delta \\n\end{array}\right)}$ (for instance Theorem 13.8 in \cite{e}).
\hspace*{\fill}$\square$\\

This leads to a genericity result for strictly under-determined matrices with respect to the Zariski topology:

\begin{thm} The behavior of a generic strictly under-determined matrix in $\mathcal{M}_{\ell < k}$ is controllable, for this set of matrices contains a Zariski open set.
\end{thm}

\noindent Proof: By Lemma 6.2, two or more generators of the cancellation ideal $\frak{i}_\ell$ of a generic strictly under-determined matrix are nonzero. For such matrices, $\dim (\mathcal{V}(\frak{i}_\ell))=n-1$ if and only if there is a nonconstant polynomial which divides all these generators. The complement of this set of matrices contains a Zariski open set by Lemma 6.3. 
\hspace*{\fill}$\square$\\

The arguments in the case of over-determined matrices are similar. Thus, now let $\ell \geq k$, and let $\mathcal{M}_{\ell, k}(d)$ again be the set of $\ell \times k$ matrices with entries in $\mathcal{A}$ of degree less than or equal to $d$. Let $N(d)= \ell k{{\tiny \left(\begin{array}{c}n+d \\n\end{array}\right)}}$ be as before, but now let $M(d)={\tiny \left(\begin{array}{c}\ell \\k \end{array}\right)}{\tiny \left(\begin{array}{c}n+kd \\n\end{array}\right)}$, so that the $k$-th determinantal ideal of an element in $\mathcal{M}_{\ell, k}(d)$ is a point in $\mathbb{C}^{M(d)}$. Thus the $k$-th determinantal ideal is a morphism $\frak{i}_k(d):\mathbb{C}^{N(d)} \rightarrow \mathbb{C}^{M(d)}$ which is algebraic, and an argument identical to the above shows:

\begin{lemma} The condition $\frak{i}_k \neq 0$ is satisfied by a Zariski open subset of $\mathcal{M}_{\ell, k}$, $\ell \geqslant k$.
\hspace*{\fill}$\square$
\end{lemma}

Finally, allowing $\ell$ to equal $k, (k+1), (k+2), ~\ldots $ ~shows that the above condition is satisfied by a generic over-determined matrix in the disjoint union $\mathcal{M}_{\ell \geqslant k} = \bigcup_{\ell \geqslant k} \mathcal{M}_{\ell,k}$. Combining this lemma with Theorem 4.2 gives:

\begin{thm} The behavior of a generic over-determined matrix in $\mathcal{M}_{\ell \geqslant k}$ is uncontrollable, namely those belonging to a Zariski open set given by the condition $\frak{i}_k \neq 0$. In particular, the behavior of a generic square matrix is uncontrollable.
\hspace*{\fill}$\square$ \\
\end{thm}

These genericity results in terms of matrices are now used to obtain genericity results for behaviors. Consider, for the sake of simplicity, behaviors defined in the space of distributions $\mathcal{D}'$ or smooth functions $\mathcal{C}^\infty$. As these are both injective cogenerators as modules over $\mathcal{A}$, there is a bijective correspondence between submodules of $\mathcal{A}^k$ and behaviors in $\mathcal{F}^k$ ($~\mathcal{F} = \mathcal{D}'$ or $\mathcal{C}^\infty$). Hence, to define a topology on the set of behaviors in these $\mathcal{F}^k$, it suffices to define it on submodules of $\mathcal{A}^k$. 

Consider first the case of under-determined behaviors. Motivated by the considerations explained in the beginning of this section, it is necessary to restrict to the subset $\mathcal{R}_{\ell,k}$ of $\mathcal{M}_{\ell,k}$ consisting of those matrices whose $\ell$ rows is a minimum set of generators for the submodules of $\mathcal{A}^k$ that they generate (Definition 2.1). Thus the submodule generated by the rows of a $P(\partial) \in \mathcal{R}_{\ell,k}$ cannot be generated by less than $\ell$ elements.
%As before, let $\mathcal{R}_{\ell,k}(d)$ be those matrices in $\mathcal{R}_{\ell,k}$ all whose entries are bounded in degree by $d$. Then the complement of $\mathcal{R}_{\ell,k}(d)$ in $\mathcal{M}_{\ell,k}(d)$, namely those $\ell \times k$ matrices (with entries bounded in degree by $d$), whose rows generate submodules that can be also generated by fewer than $\ell$ elements (entries bounded in degree by $d$), is precisely the union of the images of 
%\[
%\begin{array}{lcc}
%\mathcal{M}_{\ell, \ell -1}(r) ~\times ~ \mathcal{M}_{\ell - 1,k}(d-r) & \longrightarrow & \mathcal{M}_{\ell,k}(d) \\ \phantom{xxxx} (A(\partial) , P(\partial)) & \mapsto & A(\partial)P(\partial)
%\end{array}
%\]
%for $r=0, 1, \ldots ,d$. The direct limit of this union, as $d \rightarrow \infty$, is precisely the complement of $\mathcal{R}_{\ell,k}$ in $\mathcal{M}_{\ell,k}$. In fact, the images of 
%\[
%\mathcal{M}_{\ell,\ell - 1}(d) ~\times ~\mathcal{M}_{\ell - 1,k}(d) \longrightarrow \mathcal{M}_{\ell,k}(2d)
%\]
%is cofinal in the above direct limit, and the direct limit of this cofinal family is also the complement of $\mathcal{R}_{\ell,k}$ in $\mathcal{M}_{\ell,k}$.
Clearly, the $\ell$-th determinantal ideal $\frak{i}_\ell$ of an element in the complement of $\mathcal{R}_{\ell,k}$ in $\mathcal{M}_{\ell,k}$ equals 0, hence $\mathcal{R}_{\ell,k}$ contains a Zariski open set of $\mathcal{M}_{\ell,k}$ by Lemma 6.1.

Let $\mathcal{S}_\ell$ be the set of submodules of $\mathcal{A}^k$ which can be generated by $\ell$ ($< k$) elements, but not fewer. Consider the map 
\[
\begin{array}{lccc}
 \Pi: & \mathcal{R}_{\ell,k} & \longrightarrow & \mathcal{S}_\ell \\
  & P(\partial) & \mapsto & \mathcal{P}
 \end{array}
\]
$\ell = 1, \ldots ,k-1$, which maps a matrix in $\mathcal{R}_{\ell,k}$ to the submodule of $\mathcal{A}^k$ generated by its rows. This is a surjective map, 
and the quotient topology on $\mathcal{S}_\ell$, descending from the Zariski topology on the subspace $\mathcal{R}_{\ell,k} \subset \mathcal{M}_{\ell,k}$, is the strongest topology such that $\Pi$ is continuous. This topology on $\mathcal{S}_\ell$ is also called the Zariski topology. The disjoint union $\mathcal{S}_{\ell<k} = \bigcup_{\ell = 1}^{k-1}\mathcal{S}_\ell$ defines the set of strictly under-determined behaviors, and carries the above Zariski topology.

Theorem 6.1 now descends to behaviors:

\begin{thm} The behavior of a generic strictly under-determined submodule in $\mathcal{S}_{\ell < k}$ is controllable, for this set of behaviors contains a Zariski open set.
\hspace*{\fill}$\square$
\end{thm}

The case of over-determined behaviors is analogous. Thus, let $\mathcal{S}_\ell$ be the set of submodules of $\mathcal{A}^k$ which can be generated by $\ell$ ($\geqslant k$) elements, but not fewer. The Zariski topology on $\mathcal{S}_\ell$ is again the quotient topology it inherits from the corresponding projection $\Pi: \mathcal{R}_{\ell,k} \rightarrow \mathcal{S}_\ell$. Let $\mathcal{S}_{\ell \geqslant k}$ be the disjoint union $\bigcup_{\ell \geqslant k} \mathcal{S}_\ell $. Theorem 6.2 then implies:

\begin{thm} The behavior of a generic over-determined behavior in $\mathcal{S}_{\ell \geqslant k}$ is uncontrollable, namely those belonging to a Zariski open set given by the condition $\frak{i}_k \neq 0$. In particular, the behavior of a generic square behavior is uncontrollable.
\hspace*{\fill}$\square$ \\
\end{thm}

\section{Coordinate controllability}
All the results in this paper on the Hautus test, and on genericity, depend on conditions on the two ideals of importance, namely the cancellation ideal $\frak{i}_\ell$ and the characteristic ideal $\frak{i}_k$ of the behavior. This section introduces a weak notion of controllability suggested by these results. It also collects a few additional comments, some of which are well known to experts, but which, for the sake of completeness, are included here.

An argument similar to the proof of Theorem 3.2 (see for instance \cite{sc}), implies further that 
\[
{\frak{i}_k} \subseteq \ann(\mathcal{A}^k/\mathcal{P}) \subseteq \sqrt{\frak{i}_k}
\]
Thus $\frak{i}_k = 0$ is equivalent to $\ann(\mathcal{A}^k/\mathcal{P})= 0$.
The stronger condition, that $\mathcal{A}^k/\mathcal{P}$ be torsion free, is equivalent to controllability of $\mathcal{B}_\mathcal{F}(\mathcal{P})$, namely Theorem 1.2. What does this weaker condition imply about the behavior?

\begin{proposition} Let $\mathcal{P} \subset \mathcal{A}^k$, and $\mathcal{B}_\mathcal{F}(\mathcal{P})$ its behavior in $\mathcal{F}^k$, where $\mathcal{F}$ is either $\mathcal{D}'$, $\mathcal{C}^\infty$ or $\mathcal{S}'$. Let $\pi_j:\mathcal{F}^k \rightarrow \mathcal{F}, ~(f_1, \ldots, f_k) \mapsto f_j$ be the projection onto the $j$-th coordinate. Then $\ann(\mathcal{A}^k/\mathcal{P}) = 0$ if and only if there is a $j, ~1\leqslant j \leqslant k$, such that the restriction $\pi_j:\mathcal{B}_\mathcal{F}(\mathcal{P}) \rightarrow \mathcal{F}$ is surjective, and hence that $\pi_j(\mathcal{B}_\mathcal{F}(\mathcal{P}))$ is controllable.
\end{proposition}

\noindent Proof: Let $i_j:\mathcal{A} \rightarrow \mathcal{A}^k, ~a \mapsto (0, \ldots, a, \ldots, 0)$ be the inclusion in the $j$-th coordinate. By Proposition 4.1 in \cite{sg}, $\pi_j(\mathcal{B}_\mathcal{F}(\mathcal{P})) = \mathcal{B}_\mathcal{F}(i_j^{-1}(\mathcal{P}))$. Now, $\ann(\mathcal{A}^k/\mathcal{P}) = 0$ is equivalent to saying that there is a $j, ~1\leqslant j \leqslant k$, such that $i_j^{-1}(\mathcal{P}) = 0$. By the proposition quoted above, this is equivalent to $\pi_j(\mathcal{B}_\mathcal{F}(\mathcal{P})) = \mathcal{F}$. (In this, the fact that $\mathcal{F}$ is an injective $\mathcal{A}$-module is crucial.)
\hspace*{\fill}$\square$\\

This weaker notion may be called {\it coordinate controllability}; hence a behavior $\mathcal{B}_\mathcal{F}(\mathcal{P})$ is coordinate-controllable if and only if $\ann(\mathcal{A}^k/\mathcal{P}) = 0$. For instance, suppose that a submodule $\mathcal{P}$ of $\mathcal{A}^k$ can be generated by some $\ell$ elements, $\ell < k$. This is the case of a strictly under-determined system, and then $\frak{i}_k = 0$. Thus

\begin{corr} A strictly under-determined system is coordinate-controllable.
\hspace*{\fill}$\square$
\end{corr}

In addition, the discussion preceding Theorem 4.2 implies the following counterpart to it -- henceforth, as then, assume that $\mathcal{F}$ equals $\mathcal{D}'$ or $\mathcal{C}^\infty$.

\begin{proposition} A strictly under-determined system always contains nonzero controllable sub-behaviors, and thus it contains a unique (nonzero) maximal controllable sub-behavior.
\end{proposition}

\noindent Proof: As $\frak{i}_k = 0$, so is $\ann(\mathcal{A}^k/\mathcal{P})$, and hence 0 is an associated prime of $\mathcal{A}^k/\mathcal{P}$. \hspace*{\fill}$\square$\\

To formulate behavioral controllability (Definition 1.1) in these terms, consider a general homothety: given $r(\partial) = (a_1(\partial), \ldots ,a_k(\partial))$ in $\mathcal{A}^k$, let $i_r: \mathcal{A} \rightarrow \mathcal{A}^k$ be the morphism $a(\partial) \mapsto a(\partial)r(\partial)$ (which is an injection when $r(\partial) \neq 0$). Applying $\homo_\mathcal{A}(-,\mathcal{F})$ to this morphism gives
$r(\partial): \mathcal{F}^k \rightarrow \mathcal{F}$, which maps $f=(f_1,\ldots ,f_k)$ to $a_1(\partial)f_1 + \ldots + a_k(\partial)f_k$ (and which is a surjection when $r(\partial) \neq 0$); thus, the injections $i_j$ and the projections $\pi_j$ in Proposition 7.1 correspond to $r(\partial) = (0, \ldots, 1, \ldots, 0)$. Behavioral controllability is then an assertion about the restrictions of the maps $r(\partial)$ above to $r(\partial):\mathcal{B}_\mathcal{F}(\mathcal{P}) \rightarrow \mathcal{F}$, for all $r(\partial)$ in $\mathcal{A}^k$.

\begin{proposition} Let $\mathcal{P}$ be a submodule of $\mathcal{A}^k$. The behavior $\mathcal{B}_\mathcal{F}(\mathcal{P})$ is controllable if and only if every image $r(\partial)(\mathcal{B}_\mathcal{F}(\mathcal{P}))$ is controllable (in other words if and only if every $r(\partial)(\mathcal{B}_\mathcal{F}(\mathcal{P}))$ is either 0 or all of $\mathcal{F}$).
\end{proposition}

\noindent Proof: Consider the map 
\[\begin{array}{lcc}
\mathcal{A} & \longrightarrow & \mathcal{A}^k/\mathcal{P}  \\
 1 & \mapsto & r(\partial) + \mathcal{P} 
 \end{array}
\]
It is the 0 map if $r(\partial)$ is in $\mathcal{P}$; otherwise it is injective as $\mathcal{A}^k/\mathcal{P}$ is torsion free. Hence, applying $\homo_\mathcal{A}(-, ~\mathcal{F})$ to the above map implies the proposition.  
\hspace*{\fill}$\square$\\

A coordinate onto which a behavior surjects can play the role of an {\it input} \cite{w}. This raises the question: can the projection of a controllable behavior to several coordinates also be surjective? The example below from \cite{wooo}, shows that in general this need not be the case. \\

\noindent Example: Let $\mathcal{A} = \mathbb{C}[\partial_1,\partial_2]$. The behavior given by the kernel of the map $(\partial_1, \partial_2):\mathcal{F}^2 \rightarrow \mathcal{F}$ is controllable (for instance, its cancellation variety, the origin in $\mathbb{C}^2$, is of codimension 2). By Proposition 7.1, each $\pi_j:\mathcal{B}_\mathcal{F}((\partial_1,\partial_2))  \rightarrow \mathcal{F}, ~j = 1,2$, ~is surjective, but $\mathcal{B}_\mathcal{F}((\partial_1,\partial_2)) \rightarrow \mathcal{F}^2$ is clearly not.

A strongly controllable behavior (namely, the discussion after Proposition 2.2) does however surject onto several of its coordinates (equal in number to the rank of the free module $\mathcal{A}^k/\mathcal{P}$).
\hspace*{\fill}$\square$\\

Thus, three graded notions of controllability of a behavior $\mathcal{B}_\mathcal{F}(\mathcal{P})$ have now been identified, 
each corresponding to increasingly weaker requirements on $\mathcal{A}^k/\mathcal{P}$. 

\begin{table}[htdp]
\caption{notions of controllability}
\begin{center}
\begin{tabular}{cc}  $\phantom{--}$ strong controllability $\iff  \mathcal{A}^k/\mathcal{P}$ is free \\
$\phantom{----}$ behavioral controllability $\iff  \mathcal{A}^k/\mathcal{P}$ is torsion free\\
$\phantom{-}$coordinate controllability $\iff \ann (\mathcal{A}^k/\mathcal{P}) = 0$
\end{tabular}
\end{center}
\end{table}

Because these notions can be understood in terms of projections of the behavior to various coordinates (as the above propositions indicate), they amount to notions of inputs and outputs, and hence to a notion of causality, but will not be pursued any further here.

\vspace{.5cm}

\noindent Acknowledgement: I am grateful to Krishna Hanumanthu, Manoj Kummini, Diego Napp, Suresh Nayak and Arul Shankar for many useful conversations. I am also grateful to Maria Elena Valcher and the University of Padova for their hospitality when a part of this paper was written. Finally, I am grateful to the reviewers for their very careful reading of the manuscript.


\vspace{.5cm}
\begin{thebibliography}{99}

{\small

\bibitem{bm} H.~Bourles and B.~Marinescu, {\em Linear Time-Varying Systems}, LNCIS 410, Springer, 2011. 

\bibitem{e} D.~Eisenbud, {\em Commutative Algebra With A View Toward Algebraic Geometry}, Graduate Texts in Mathematics 150, Springer-Verlag, 1995. 

\bibitem{g} H.~Gl\"using-L\"uerssen, {\em Linear Delay-Differential Systems with Commensurate Delays: An Algebraic Approach}, Lecture Notes in Mathematics 1770, Springer 2002.

\bibitem{h} M.L.J.~Hautus, {\em Controllability and observability conditions for linear autonomous systems}, Nederlandse Akademie van Wetenschappen, Proc. Ser. A, 72:443-448, 1969.

\bibitem{k} R.E.~Kalman, {\em On the general theory of control systems}, 
Proceedings of the 1st World Congress of the International Federation of
Automatic Control, 481-493, Moscow, 1960.

\bibitem{l} V.~Lomadze, {\em Controllability as minimality}, SIAM jl. Control and Optimization, 50:357-367, 2012. 

\bibitem{l1} V.~Lomadze, {\em The PBH test for LTID systems in several variables}, preprint, 2013.

%\bibitem{m} H.~Matsumura, {\em Commutative Ring Theory}, Cambridge Studies
%in Advanced Mathematics, 1989.

\bibitem{nps} D.~Napp, M.~van der Put and S.~Shankar, {\em Periodic behaviors}, SIAM jl. Control and Optimization, 48:4652-4663, 2010.

\bibitem{o} U.~Oberst, {\em Multidimensional constant linear systems}, 
Acta Applicandae Mathematicae, 20:1-175, 1990.  

\bibitem{ps} H.K.~Pillai and S.~Shankar, {\em A behavioural approach to control of distributed systems}, SIAM jl. Control and Optimization, 37:388-408, 1998.

\bibitem{pq} J-F.~Pommaret and A.~Quadrat, {\em Algebraic analysis of 
linear multidimensional control systems}, IMA jl. on 
Mathematical Control and Information, 16:275-297, 1999.

\bibitem{rj} P.~Rocha and J.C.~Willems, {\em Behavioral controllability of delay-differential systems}, SIAM jl. Control and Optimization, 35:254-264, 1997.

\bibitem{rw} P.~Rocha and J.~Wood, {\em Trajectory control and interconnection of 1D and nD systems}, SIAM jl. Control and Optimization, 40:107-134, 2001.

%\bibitem{sn} S.~Shankar, The Nullstellensatz for systems of PDE, {\em Advances in Applied Mathematics}, 23:360-374, 1999. 

\bibitem{se} S.~Shankar, {\em The evolution of the concept of controllability}, Mathematical and Computer Modelling of Dynamical Systems, vol.8, 397-406, 2002.

\bibitem{sg} S.~Shankar, {\em Geometric completeness of distribution spaces}, Acta Applicandae Mathematicae, 77:163-180, 2003.

\bibitem{sc} S.~Shankar, {\em A Cousin problem for systems of PDE}, Mathematische Nachrichten, 280:446-450, 2007.

\bibitem{w} J.C.~Willems, {\em The behavioral approach to open and interconnected systems}, IEEE Control Systems Magazine, 27:46-99, 2007.

\bibitem{wo} W.A.~Wolovich, {\em Linear Multivariable Systems}, Springer, 1974.

\bibitem{woo} J.~Wood, E.~Rogers and D.H.~Owens, {\em A formal theory of matrix primeness}, Mathematics of Control, Signals and Systems, 11:40-78, 1998.

\bibitem{wooo} J.~Wood, U.~Oberst, E.~Rogers and D.H.~Owens, {\em A behavioral approach to the pole structure of one-dimensional and multidimensional linear systems}, SIAM jl. Control and Optimization, 38:627-661, 2000.

\bibitem{z} E.~Zerz, {\em Primeness of multivariate polynomial matrices}, Systems Control Letters, 29:139-145, 1996.
}
\end{thebibliography}
\end{document}